\documentclass{elsart}
\usepackage{epsfig}
\usepackage{amsfonts}
\usepackage{amssymb}

\def\asin{\hbox{\rm asin}}

\def\Card{{\hbox{\rm  Card\,}}}

\def\ord{\hbox{\rm ord}}
\def\rank{\hbox{\rm Rank}}
\def\inter{\hbox{\rm int}}
\def\bZd{\Bbb {Z}^d}
\def\bZ{\Bbb Z}

\def\cal{\mathbf}
\newcommand{\nn}{\nonumber \\}
\begin{document}
\begin{frontmatter}
\title{Nonparametric Denoising of Signals with Unknown Local Structure, I: Oracle Inequalities}
\author[Tolya]{Anatoli Juditsky\corauthref{cor}},
\corauth[cor]{Corresponding author.}
\ead{anatoli.juditsky@imag.fr}
\author[Arik]{Arkadi Nemirovski}
\ead{nemirovs@isye.gatech.edu}
\address[Tolya]{LJK, B.P. 53, 38041 Grenoble Cedex 9, France}
\address[Arik]{ISyE,
Georgia Institute of Technology,
765 Ferst Drive,
 Atlanta GA 30332-0205 USA}
\begin{abstract}
We consider the problem of pointwise estimation of multi-dimensional signals
$s$,
from noisy observations $(y_\tau)$ on the regular grid $\bZd$.
Our focus is on the adaptive estimation in the case
when the signal can be well recovered using a (hypothetical)
linear filter, which can depend on the unknown
signal itself.
\par
The basic setting of the problem we address here can be summarized as follows:
suppose that the signal $s$ is ``well-filtered'', i.e.
there exists an adapted time-invariant linear filter $q^*_T$ with the coefficients
which vanish outside
the ``cube'' $\{0,..., T\}^d$ which recovers $s_0$ from observations
with small mean-squared error.
We suppose that we do not know the filter $q^*$, although, we do know that such a filter exists.
We give partial answers to the following questions:
\begin{description}
\item[]-- is it possible to construct an adaptive estimator of the value $s_0$, which relies upon observations
and recovers $s_0$ with basically the same estimation error as the unknown filter $q^*_T$?
\item[]-- how rich is the family of well-filtered (in the above sense) signals?
\end{description}
We show that the answer to the first question is affirmative and provide a numerically efficient
construction of a nonlinear adaptive filter. Further, we establish a simple calculus of ``well-filtered" signals,
and show that their family is quite large: it contains, for instance,
sampled smooth signals, sampled modulated smooth signals and sampled harmonic functions.
\end{abstract}
\end{frontmatter}

\section{Introduction}
 In this paper, we focus on the problem of denoising of multi-dimensional signals.
 Let ${\cal F}=(\Omega,\Sigma,P)$ be a probability space. We
consider the problem of recovering unknown random field
$(s_\tau=s_\tau(\xi))_{{\tau\in\bZd\atop\xi\in\Omega}}$ over
$\bZd$ from noisy observations
\begin{equation}\label{eq:mod}
y_{\tau}=s_{\tau}+e_{\tau}.
\end{equation}
 It is
convenient for us to assume that both the {\sl signal} $(s_\tau)$
and the noises are complex-valued. Besides this, we assume that
the field  $(e_{\tau})$ of observation noises is independent of
$(s_{\tau})$ and is of the form $e_{\tau}=\sigma \epsilon_{\tau}$,
where $(\epsilon_{\tau})$ are independent of each other {\sl
standard} Gaussian complex-valued variables; the adjective
``standard'' means that $\Re(\epsilon_{\tau})$,
$\Im(\epsilon_{\tau})$ are independent of each other ${\cal
N}(0,1)$ random variables. Our focus here is at estimating the value $s_t$ of signal at a given location $t\in \bZd$.
\par
The above problem is  ``classical''  in statistical estimation and signal processing, and as such,
has received much attention. In particular, linear estimators (referred as linear filters in the signal processing community) are widely used in the statistical literature. 
To be more precise, suppose that our aim is to recover the value $s_0$ of the signal at zero given  observations $(y_\tau)$ on the box ${\cal O}_T=\{\tau\in\bZd:|\tau_j|\leq T,1\leq j\leq d\}$. We call the estimation $\hat{s}$  of $s_0$  linear if it is of the form
\[
\hat{s}_\ell=\sum_{\tau\in \cal O_T} q_{\tau} y_{\tau}
\]
for some $q\in C(\cal O_T)$, where $C(\cal O_T)$ is the set of complex-valued fields $q=\{q_\tau,\;\tau\in \cal O_T\}$ over $\cal O_T$.
\par
The simplicity of linear estimators is responsible for their popularity in statistical signal processing. Another outstanding feature of such estimators is their minimax property.
 Suppose that the {\em a priori} information resumes to the fact that $(s_\tau)$ belongs to some convex compact set which is symmetric with respect to zero, let us call it $\cal S$.
One of the most renown results of estimation theory (see, for instance, \cite{IKh84,DonLow92,Don95}) states that the  {\em linear minimax estimator}
is, in a certain sense, an optimal estimator of $s_0$ in our problem. Indeed,
consider the following {\em linear minimax estimation strategy}:
 let $q^{(T)}_*$ be the optimal solution\footnote{For evident reasons such a solution exists in the situation we are interested in.} to the problem
\begin{eqnarray*}
\min_{q\in C(\cal O_T)}\; \max_{s\in \cal S} E_{s}\left(s_0-\sum_{\tau\in \cal O_T} q_\tau y_{\tau}\right)^2
\end{eqnarray*}
(here $E_s$ stands for the expectation with respect to the distribution of $(y_\tau)$ which corresponds to the underlying signal $s$). The {  linear minimax estimator} $\hat{s}^*_\ell$ of $s_0$ is defined by
\[
\hat{s}^*_\ell=\sum_{\tau\in \cal O_\tau} q^{(T)}_{*,\tau} y_{\tau}.
\]
Then
\[
\max_{s\in \cal S} E_{s}(s_0-\hat{s}^*_\ell)^2\le C \inf_{\hat{s}}\max_{s\in \cal S} E_{s}(s_0-\hat{s})^2,
\]
where the infimum in the right-hand side is taken over all possible estimators of $s_0$ from observations $(y_\tau)$ and $C$ is a moderate absolute constant (e.g., $C\le 1.25$). In other words, the linear estimator $\hat{s}_\ell$ is a (almost) minimax estimator of $s_0$. We would like to stress the exceptional power of the above result -- we only need $\cal S$ to be  convex and compact for the linear estimator to be minimax optimal. The evident downside of using linear minimax estimators is that the {\em a priori} information about the set $\cal S$ of signals should be as precise as possible to achieve descent estimation accuracy. There was a significant research on {\sl adaptive
estimation} in the above setting (cf \cite{cailow1,cailow2}). Those techniques allow to choose the ``best'' in a certains sense set which contains the signal from special finite families of convex sets.
Another ``classical'' approach to adaptation for linear estimators has been developed in \cite{{bib:lep90},bib:lep91,{bib:lep92},lep-spok97,tsyb}. In the latter approach the ``form'' of the filter $q^{(T)}$ is considered as given in advance (no information about sets of signals is used in this case), and the parameter $T$ (the ``window width'') is selected adaptively to achieve the best bias/variance tradeoff. Recently, more general adaptation techniques has been studied in \cite{LL,GL}, which allow to choose the best  estimator from special finite families of available linear estimators.
\par
The problem we are interested in here, {\sl when posed informally},  is as follows:
if we consider the form of the filter as a ``free parameter", is it possible to provide an estimation procedure
which is adaptive with respect to this parameter? In other words, suppose that a ``good" filter
$q^{(T)}_*$, with a small estimation error exists. Then, is it possible  to construct
a data-driven estimation method which has (almost) the same  accuracy as the ``oracle" -- a hypothetic optimal
estimation method which uses the ``good" filter $q^{(T)}_*$. It is natural, as it is common in adaptive nonparametric estimation,
to measure the quality of an adaptive estimation routine with the factor by which
the risk of the adaptive procedure is greater than that of the ``oracle" estimator. What we look for
is the estimation method for which this factor is not too large.
Let us consider, for instance, the following question:
\par
{\sl (?)} {\sl Suppose that know that the (deterministic or random) signal
$(s_\tau)_{\tau\in\bZd}\equiv(s_\tau(\xi))_{{\tau\in\bZd\atop\xi\in\Omega}}$
underlying observations {\rm (\ref{eq:mod})} can be recovered from
these observations ``at a parametric rate'' by ``linear
time-invariant filtering'': for a given $T$, there exists (unknown
in advance) filter $q^{(T)}_*$ which recovers
$s_0$ via $O(T^d)$ observations around zero such
that
\begin{equation}\label{asnem101}
E\left\{|s_0-\sum\limits_{\tau\in{\cal O}_T} q^{(T)}_{*,\tau}
y_{\tau}|^2\right\}\leq
O(\sigma^2T^{-d}).
\end{equation}
Can we mimic this filter?}
\par
We show that the answer to the question {\sl (?)} is positive. Namely,
whenever a discrete time signal (that is, a signal defined on
a regular discrete grid) is {\sl well-filtered}, i.e., can be
recovered from its noisy observations {\sl at a parametric rate}
by a {\sl linear time-invariant filter}, we can recover this
signal at a ``nearly parametric'' rate {\sl without a priori
knowledge of the associated filter}.
\par
Several points should be stressed in the above claim. First, we are able to mimic
only ideal filters $q^{(T)}_*$ of small  $l_2$-norm.
Indeed, the relation (\ref{asnem101}) implies that the stochastic term of the error
$E\left(\sum\limits_{\tau\in{\cal O}_T} q^{(T)}_{*,\tau}
e_{\tau}\right)^2$ is bounded with $O(\sigma^2T^{-d})$, which is conceivable only if $|q^{(T)}_*|_2=O(T^{-d/2})$.
This constraint is crucial, as the price for adaptation
becomes prohibitive when the $l_2$-norm of the ideal filter is much larger than $O(T^{-d/2})$.
Though this assumption seems quite restrictive, the family of well-filtered signals is quite wide. As we shall see later, this family contains also ``highly
oscillating" sampled modulated smooth signals, sampled harmonic functions, etc.
\par
In this paper we also treat the problem of adaptive prediction, when we are interested in
recovering of a discrete time signal at a point $t\in\bZd$ via
noisy observations taken at the points $\{\tau\in\bZd:t_j-T\leq
\tau_j\leq t_j-\kappa\}$ ``preceding'' the point $t$, with a given
in advance ``forecast horizon''  $\kappa\geq0$.
\par
The rest of our paper is organized as follows.
In Section \ref{problem_st} we give a formal
definition of a well-filtered (well-predicted) signal on a $d$-dimensional regular
grid (the latter, w.l.o.g., is normalized to be $\bZd$), and then show
in Section \ref{main_result}
demonstrate that such a signal can be recovered at a nearly
parametric rate {\sl without a priori knowledge of the
corresponding ``good filter''} (Theorems \ref{the:Oadapt} and \ref{the:Opredict}). The
underlying estimation routines
(i.e., ``Algorithm A'' of Section \ref{AlgorithmA} and ``Algorithm B'' of Section \ref{AlgorithmB}) constitute
a substantial extension of the procedures proposed in \cite{Nem81} and  \cite{Medallion}. In Section
\ref{calculus}, we demonstrate that the family
 of well-filtered signals is
pretty wide -- it contains a wide spectrum of  ``basic functions''
(for example, exponential polynomials) and is closed with respect
to a number of basic operations, including modulation, taking
linear combinations and tensor products.
\par To make the exposition more readable,
 all proofs are collected in the appendix.
\par
The denoising procedures, described in this paper constitute the basic bricks
of the construction of adaptive estimators of {\sl locally well-filtered} signals,
which we describe in the companion paper\cite{denoising2}.
The results of \cite{denoising2}
extend to the wide classes of modulated signals
the results of
\cite{bib:pief84,bib:donjohnker93,bib:jud94,GN1}
on spatial adaptive estimates of  signals with inhomogeneous smoothness.
%
%
\section{Problem statement}\label{problem_st}
In order to proceed we need some notations.\\
{\sl Fields over $\bZd$.} Let $ C({\Bbb  Z}^d)$ be the linear
space of complex-valued fields $r=\{r_{\tau}:\tau\in\bZd\}$
over $\bZd$. \\
$\bullet$
Given nonnegative integer $T$ and $p\in[1,\infty]$, we
define semi-norms $|\cdot|_{T,p}$ on $ C({\Bbb  Z}^d)$ by $
|r|_{T,p} = \left(\sum\limits_{|\tau|\leq T}
|r_{\tau}|^p\right)^{1/p}$,
$|\tau|=\max\{|\tau_1|,...,|\tau_d|\}$, with the standard
interpretation of the right hand side when $p=\infty$, and we set
$ |r|_p =\lim_{T\to\infty} |r|_{T,p}\in {\Bbb R}\cup\{+\infty\}. $
A field $r\in  C({\Bbb  Z}^d)$ with finitely many nonzero entries
$r_{\tau}$ is called {\sl a filter}, and the smallest $T$ such
that $r_{\tau}=0$ whenever $|\tau|>T$, is called the {\sl order}
$\ord(r)$ of a filter $r$; we write $ C_T(\bZd)=\{r\in C({\Bbb
Z}^d)\mid\, \ord(r)\leq T\}. $
 We identify a filter
$r$ with the multivariate Laurent sum $
r(z_1,...,z_d)=\sum\limits_{\tau}
r_{\tau}z_1^{\tau_1}...z_d^{\tau_d}. $\\ $\bullet$ We call a
filter $r$ {\sl polynomial}, if the corresponding Laurent sum is a
polynomial (i.e., if the  entries $r_{\tau}$ vanish when any of $\tau_j< 0$, $j=1,...,d$). The set of all polynomials is
denoted $P(\bZd)$. For integers $k,T$, $0\leq k\leq T$, we denote
by $P_T^k(\bZd)$ the subspace of $P(\bZd)$ formed by polynomials
$r$ for which the entries $r_\tau$ vanish outside the set
$k\le \tau_j\le T$, $j=1,...,d$.
\\
$\bullet$ We denote by $\Delta_j$, $j=1,...,d$, the ``basic shift
operators'' on $C({\Bbb  Z}^d)$:
$$
(\Delta_jr)_{\tau_1,...,\tau_d}=r_{\tau_1,...,\tau_{j-1},\tau_j-1,\tau_{j+1},...,\tau_d}.
$$ Further, we use the notation $\Delta_j^{-1}$ for the inverse of $\Delta_j$:
$$
(\Delta_j^{-1}r)_{\tau_1,...,\tau_d}=r_{\tau_1,...,\tau_{j-1},\tau_j+1,\tau_{j+1},...,\tau_d}.
$$
$\bullet$ Finally, we define the output of a filter  $r$, the
input to the filter being a field $x\in C(\bZd)$, as the field $
r(\Delta)x\equiv r(\Delta_1,\Delta_2,...,\Delta_d)x, $ so that $
(r(\Delta)x)_t=\sum\limits_\tau r_\tau x_{t-\tau}. $
\\
{\sl Fourier transform.} Let $T$ be a nonnegative integer, let
$\Gamma_T$ be the set of roots of 1 of the degree $2T+1$, and let
$C(\Gamma_T^d)$ be the space of complex-valued functions on
$\Gamma_T^d\equiv(\Gamma_T)^d$.
\\
$\bullet$ We define the Fourier transform $F_T: C({\Bbb  Z}^d)\to
C(\Gamma_T^d)$ as $(F_Tr)(\mu)= {1\over (2T+1)^{d/2}}
\sum\limits_{|\tau|\leq T}
r_{\tau}\mu_1^{\tau_1}...\mu_d^{\tau_d}\equiv{1\over (2T+1)^{d/2}}
r(\mu), r\in C_T(\bZd)$, where $\mu\in\Gamma_T^d$. Note that $
r_{\tau} = {1\over(2T+1)^{d/2}}\sum\limits_{\mu\in\Gamma_T^d}
(F_Tr)(\mu)\mu_1^{-\tau_1}...\mu_d^{-\tau_d},\,\,\forall
(\tau:|\tau|\leq T). $ The Fourier transform allows to equip $
C({\Bbb  Z}^d)$ with semi-norms coming from the standard $p$-norms
on $C(\Gamma_T^d)$:
$$ |r|_{T,p}^* = |F_Tr|_p
\equiv
\left(\sum\limits_{\mu\in\Gamma_T^d}|(F_Tr)(\mu)|^p\right)^{1/p},
$$ with the standard interpretation of the right hand side for
$p=\infty$.
\par
Now it is time to give a precise meaning to the basic question {\sl (?)} of Introduction.
In order to do this, we should specify our a priori
knowledge of the constant factor hidden in $O(\cdot)$ and on the
ranges on values of $T$ and $\tau$ where (\ref{asnem101}) holds
true.
\subsection{Nice signals}
 Since the observation noises are independent of $(s_\tau)$,
we have
\begin{equation}\label{asnem102}
E\left\{|s_\tau-(q(\Delta)y)_\tau|^2\right\}=2\sigma^2|q|_2^2+
E_\xi\left\{|s_\tau(\xi)-(q(\Delta)s(\xi))_\tau|^2\right\};
\end{equation}
therefore in order to ensure (\ref{asnem101}), both terms in the
right hand side of the latter inequality should be of order of
$T^{-d}$. This observation motivates the following
\begin{defn}\label{def:wfs} Let $\theta\geq0$,
$\rho\geq1$ be reals, let $L$ be a nonnegative integer or
$+\infty$, and let $t\in\bZd$. Finally, let
$(s_\tau)_{\tau\in\bZd}\equiv(s_\tau(\xi))_{{\tau\in\bZd\atop\xi\in\Omega}}$
be a random field on $\bZd$.
\\
(1) {\rm [$T$-well-filtered signals]} Let $T$ be a nonnegative
integer. We say that $(s_\tau)$ is $T$-well-filtered, with the
parameters $\theta$, $\rho$, $L$, at the point $t$ {\rm (notation:
$(s_\tau)\in{\cal S}^{t}_L(\theta,\rho,T)$)}, if there exists a
filter $q=q^{(T)}\in C_T(\bZd)$, $|q|_2\le {\rho\over
(2T+1)^{d/2}}$, which reproduces $(s_\tau)$ in the box
$\{\tau:|\tau-t|\leq L\}$ with the mean square error not exceeding
$\theta(2T+1)^{-d/2}$:
\begin{equation}\label{errorsmall}
\max\limits_{\tau:|\tau-t|\leq L}\left[E\left\{\left|s_{\tau}-
(q(\Delta)s)_{\tau}\right|^2\right\}\right]^{1/2}\leq
\theta(2T+1)^{-d/2}.
\end{equation}
(2) {\rm [well-filtered signals]} We say that $(s_\tau)$ is
well-filtered,  with the parameters $\theta$, $\rho$, $L$, at the
point $t$ {\rm (we use the notation: $(s_\tau)\in{\cal
F}^t_L(\theta,\rho)$)}, if, for every integer $T$, $0\leq T\leq
L$, $(s_\tau)$ is $T$-well-filtered, with the parameters
$\theta,\rho,L$, at $t$.
\end{defn}
 In the above definition we were focusing on the case of {\sl
de-noising} -- recovering a well-filtered signal $(s)$ at a point
$t\in\bZd$ via a given number observations ``around'' this point.\footnote{To be more precise,
in the filtering literature  this case
is referred to as {\em interpolation}.}
Another interesting problem is that of {\sl prediction}, where the
goal is to recover $s_t$ via observations $y_\tau$ ``preceding by
a given horizon $\kappa\in \bZ_+$'' the point $t$, i.e.,
observations with $\tau_j\leq t_j-\kappa$, $j=1,...,d$.
\begin{defn}\label{defpredict}
 Let $\theta\geq0$,
$\rho\geq1$ be reals, let $T_0\geq\kappa$ be nonnegative integers,
$L$ be a nonnegative integer or $+\infty$, and let $t\in\bZd$.
Finally, let
$(s_\tau)_{\tau\in\bZd}\equiv(s_\tau(\xi))_{{\tau\in\bZd\atop\xi\in\Omega}}$
be a random field on $\bZd$.
\\
(1)
 {\rm [$T$-well-predicted signals]}
Let $T$ be a nonnegative integer. We say that $(s_\tau)$ is
$T$-well predicted with the parameters $\theta$, $\rho$, $\kappa$,
$L$, at the point $t$ {\rm (notation: $(s_\tau)\in{\cal
Q}^{t}_{\kappa,L}(\theta,\rho,T)$)}, if there exists a filter
$q=q^{(T)}\in P_T^\kappa(\bZd)$, $|q|_2\leq {\rho\over
(2T+1)^{d/2}}$, which reproduces $(s_\tau)$ in the box
$\{\tau:|\tau-t|\leq L\}$ with the mean square error not exceeding
$\theta(2T+1)^{-d/2}$:
\begin{equation}
\label{errorsmallpred} \max\limits_{\tau:|\tau-t|\leq
L}\left[E\left\{\left|s_{\tau}-
(q(\Delta)s)_{\tau}\right|^2\right\}\right]^{1/2}\leq
\theta(2T+1)^{-d/2}.
\end{equation}
(2)
 {\rm [well-predicted signals]} We say that $(s_\tau)$ is
well-predicted,  with the parameters $\theta$, $\rho$, $\kappa$,
$T_0$, $L$, at the point $t$ {\rm (notation: $(s_\tau)\in{\cal
P}^t_{\kappa,T_0,L}(\theta,\rho)$)}, if, for every integer $T$,
$T_0\leq T\leq L$, $(s_\tau)$ is $T$-well-predicted, with the
parameters $\theta,\rho,\kappa,L$, at $t$.
\end{defn}
\begin{rem}\label{rem:wellpred} Note that the quantitative
description of a well-predicted field, when compared with the
description of a well-filtered field, involves an extra parameter
$T_0$ -- the smallest ``window width'' starting with which a
possibility to predict $s_t$ is postulated. In the case of
well-filtered fields, this width is just 0, in full accordance
with the fact that in the de-noising problem every signal is
$0$-well-filtered, at every point $t$, with parameters $\theta=0$,
$\rho=1$, $L=\infty$ due to the existence of the trivial
``single-point'' filter $q(z)\equiv 1$. \end{rem}

In the sequel, we qualify as {\sl nice} a signal
which fulfils the requirements of Definition \ref{defpredict} or \ref{def:wfs} above.
The filters $q^{(T)}$ associated, in the
sense of the above definitions, with a nice signal
$(s_\tau)$ as to filters {\sl certifying} the
``niceness"  (``well-filterability" of ``well-predictability") of the signal.
\\
 We are about to demonstrate that in the framework,
suggested by the above definitions,
the answer to the question {\sl (?)} is affirmative.I.e.,  a signal which is
nice ($T$-well-filtered or $T$-well-predicted, with parameters $\theta,\rho,L=3T$) at a point
$t$ can be recovered at this point ``at a nearly parametric
rate'' with {\sl no a priori knowledge of the corresponding ``good
filter''}; all we should know in advance are the parameters $\rho$ and $T$.
\section{Main result}
\label{main_result}
We
start the recovering routine for the adaptive filtering problem.
\subsection{Adaptive filtering}
\label{AlgorithmA} The estimator we intend to use is as follows:
\\
\textbf{Algorithm A:} {\sl Given a setup $(\rho\geq1,T)$ and a
point $t\in\bZd$, we build an estimation $\widehat{s}_{t}[T,y]$ of
$s_{t}$ via observations $(y_{\tau})$, $|\tau-t|\leq 4T$,  as
follows:
\\
(1) When $T=0$, we merely set $\widehat{s}_{t}[0,y]=y_t$
\\
(2) When $T>0$, we set
$\widehat{s}_{t}[T,y]=(\widehat{\phi}^{t}(\Delta)y)_{t},$ where
$\widehat{\phi}^{t}\in C_{2T}(\bZ^2)$ is an optimal solution to
the following optimization problem: }
\begin{equation}
\min\limits_{\phi\in
C_{2T}(\bZd)}\bigg\{\underbrace{|\Delta_1^{-t_1}...\Delta_d^{-t_d}(1-\phi(\Delta))y|^*_{2T,\infty}}_{
J(\phi,y^{t}_{4T})}: |\phi|^*_{2T,1}\leq
2^{d/2}\rho^2(2T+1)^{-d/2} \bigg\},\;
\label{eq:optim}
\end{equation}
{\sl where}
$
y^{t}_L=\left\{y_{\tau}: |t-\tau|\leq L\right\}.
$
\par
 Note that the objective in (\ref{eq:optim}) is affected
only by observations  $y^{t}_{4T}$, so that our algorithm recovers
$s_{t}$ via $(8T+1)^d$ observations ``around'' the point $t$.
\begin{thm}
\label{the:Oadapt} Assume that the signal $(s_\tau)$ underlying
observations {\rm (\ref{eq:mod})} is $T$-well-filtered, with
parameters $\theta$, $\rho$, $L\geq 3T$:  $(s_\tau)\in{\cal
S}^t_L(\theta,\rho,T)$ with $L\geq 3T$. Then the mean square error
of the estimate $\widehat{s}_{t}[T,\cdot]$ of $s_{t}$ yielded by
Algorithm A with setup $(\rho,T)$ can be bounded from above as
follows:
\begin{equation}\label{eq:0th}
\begin{array}{rcl}
\left(E\left\{\left|\widehat{s}_{t}[T,y]-s_{t}\right|^2\right\}\right)^{1/2}
&\leq&\displaystyle{
c(d)\rho^3{\theta+\sigma\rho\sqrt{\ln(2T+1)+1}\over(2T+1)^{d/2}}},\\[4mm]
c(d)&=&3(2^d+2^{3d-1}).\\
\end{array}
\end{equation}
In particular, if $(s_\tau)$ is well-filtered, with the parameters
$\theta$, $\rho$, $L$, at a point $t$, then for every integer $T$,
$0\leq T\leq \lfloor L/3\rfloor$, the accuracy of the estimate
$\widehat{s}_t[T,y]$ of $s_t$ yielded by Algorithm A can be
bounded by {\rm (\ref{eq:0th})}. Finally, in the case of
deterministic $(s)$, we have
\begin{equation}\label{ThatIsIt}
\begin{array}{l}
|s_t-\widehat{s}_{t}[T,y]|\leq c(d)\rho^3\left[\theta+
\sigma\rho\Theta^t_T\right](2T+1)^{-d/2},\\
\Theta^t_T=\sigma^{-1}\max\limits_{\tau:|\tau|\leq 2T}
|\Delta_1^{\tau_1-t_1}...\Delta_d^{\tau_d-t_d}e|_{2T,\infty}^*.\\
\end{array}
\end{equation}
\end{thm}
{\sl Comments:} note that Theorem \ref{the:Oadapt}  gives an
affirmative answer to the question {\sl (?)}. Indeed, let a signal
$(s_\tau)$ admit, for some $T$, a filter-type estimate $
\bar{s}_\tau=(q^*(\Delta)y)_\tau$ with ``window width'' $T$ (i.e.,
with $q^*\in C_T(\bZd)$) and with the mean square error which, in
an $O(T)$-neighborhood of a point $t$, is of the ``parametric''
order $O\left(\sigma(2T+1)^{-d/2}\right)$:
\begin{equation}\label{asnem103}\max_{\tau:|\tau-t|\leq 3T}
E\left\{\left|s_\tau-\bar{s}_\tau\right|^2\right\}\leq
\kappa^2\equiv {\sigma^2\mu^2\over (2T+1)^{d/2}}\end{equation}
with some known $\mu\geq1$. We do not know what is this estimate,
although do know that it exists (i.e., know the associated
$T,\mu$), and we want to recover $s_t$ from observations
$y^t_{4T}$ nearly as well as if we were using our hypothetic
estimate $\bar{s}_t$. Theorem \ref{the:Oadapt} says that Algorithm
A  basically achieves this goal. Indeed, from (\ref{asnem102}),
(\ref{asnem103}) it follows that
  $|q^*|_2\leq {\mu\over(2T+1)^{d/2}}$ and
$(s_\tau)\in{\cal S}^t_{3T}(\sigma\mu,\mu,T)$. Applying Theorem
\ref{the:Oadapt} with $\rho=\mu$, $\theta=\sigma\mu$, $L=3T$, we
conclude that with the estimate yielded by Algorithm A, the mean
square error of recovering $s_t$ does not exceed $
O(1)\mu^3\left[1+\sqrt{\ln(2T+1)}\right]\kappa. $ We see that {\sl
as far as the dependence on ``observation time'' $T^d$ is
concerned}, the estimate yielded by Algorithm A is just by a
logarithmic in $T$ factor worse than the estimate $\bar{s}_t$ we
wish to mimic.
\\
In the literature on nonparametric estimation the bounds as in
Theorem \ref{the:Oadapt} are often referred to as {\em oracle
inequalities}. Since the pioneering work \cite{Akaike} a number of
oracle inequalities have been established for a wide variety of
estimation problems (cf. the papers \cite{Kneip},
\cite{BarBirgeMas}, \cite{BirgeLep}, \cite{BirgeMassartGauss},
\cite{bib:donjohn92b}, \cite{DonJohnKerPic}, \cite{Caiwave} among
many others). In that context one refer to the filter $q$, which
certifies the niceness of the signal, as the oracle, and
the bound (\ref{eq:0th}) describes the ability of a particular
adaptive method (Algorithm A above) to reproduce the oracle.
\par
Note that the ``upper bound'' of Theorem \ref{the:Oadapt} may be compared to the lower bound of Theorem 2 of \cite{AnnIHP} for the $1$-dimensional situation. The latter result states that  one can exhibit a family of signals which 1) each member of the family can be recovered with the rate $O\big({\sigma\rho\over \sqrt{T}}\big)$ using the corresponding certifying filter; 2) the rate of estimation of signals from the family using the observation (\ref{eq:mod}) is at best $O\Big({\sigma\rho^2 }\sqrt{\ln T\over T}\Big)$. In other words, it states that the factor $\rho
\sqrt{\ln (2T+1)}$ is an unavoidable ``price" for adaptation. When
comparing the result of Theorem \ref{the:Oadapt} to that lower bound,
we observe an extra factor $\rho^2\ge 1$
in the corresponding upper bound (\ref{eq:0th}).  By now we do not
know if this extra factor can be completely eliminated.
Nevertheless, in light of these results, we can claim that
recovering of signals with certifying filter of large $l_2$-norm
is a rather desperate task -- the price for adaptation is then
proportional to $\rho\gg 1$ in this case.
\subsection{Adaptive prediction}
\label{AlgorithmB}
\label{prediction} We now turn to the problem of adaptive prediction.
The predictor we intend to use is as
follows:
\\
\textbf{Algorithm B:} {\sl Given a setup $(\rho\geq1,\kappa,T)$
and a point $t\in\bZd$, we build a prediction
$\widehat{s}_{t}[\kappa,T,y]$ of $s_t$ via observations
$(y_{\tau})$, $\kappa\leq t_j-\tau_j\leq 4T$, $j=1,...,d$,  as
$\widehat{s}_{t}[\kappa,T,y]=(\widehat{\psi}^{t}(\Delta)y)_{t},$
where $\widehat{\psi}^{t}\in P_{2T}^\kappa(\bZ^2)$ is an optimal
solution to the following optimization problem:}
\begin{equation}\label{eq:optimpred}
\min\limits_{\psi\in C_{2T}^\kappa(\bZd)}\bigg\{
\underbrace{|\Delta_1^{-t_1}...\Delta_d^{-t_d}(1-\psi(\Delta))y|^*_{2T,\infty}}_{J(\psi,y^{t}_{\kappa,4T})}:
|\psi|^*_{2T,1}\leq {2^{d/2}\rho^2\over(2T+1)^{d/2}} \bigg\};
\end{equation}
{\sl where} $ y^{t}_{\kappa,L}=\left\{y_{\tau}: \kappa\leq
t_j-\tau_j\leq L,\,j=1,...,d \right\}. $\\
 Note that the objective in (\ref{eq:optimpred}) is affected
only by observations  $y^{t}_{\kappa,4T}$, so that our algorithm
recovers $s_{t}$ via $(4T-\kappa+1)^d$ observations ``around'' the
point $t$.
\begin{thm}
\label{the:Opredict} Assume that the signal $(s_\tau)$ underlying
observations {\rm (\ref{eq:mod})} is $T$-well-predicted, with
parameters $\theta$, $\rho$, $\kappa$, $L\geq 3T$:
$(s_\tau)\in{\cal Q}^t_{\kappa,L}(\theta,\rho,T)$ with $L\geq 3T$.
Then the mean square error of the estimate
$\widehat{s}_{t}[\kappa,T,\cdot]$ of $s_{t}$, provided by Algorithm
B with setup $(\rho,\kappa,T)$, can be bounded from above as
follows:
\begin{equation}\label{eq:0thpred}
\begin{array}{rcl}
\left(E\left\{\left|\widehat{s}_{t}[\kappa,T,y]-s_{t}\right|^2\right\}\right)^{1/2}
&\leq&\displaystyle{
c(d)\rho^3{\theta+\sigma\rho\sqrt{\ln(2T+1)+1}\over(2T+1)^{d/2}}},\\[4mm]
c(d)&=&3(2^d+2^{3d-1}).\\
\end{array}
\end{equation}
In particular, if $(s_\tau)$ is well-predicted, with the
parameters $\theta$, $\rho$, $\kappa$, $T_0$, $L$, at a point $t$,
then for every integer $T$, $T_0\leq T\leq \lfloor L/3\rfloor$,
the accuracy of the estimate $\widehat{s}_t[\kappa,T,y]$ of $s_t$
yielded by Algorithm B can be bounded by {\rm
(\ref{eq:0thpred})}.\\
Finally, in the case of deterministic $(s)$, we have
\begin{equation}\label{ThatIsItP}
\begin{array}{rcl}
|s_t-\widehat{s}_{t}[T,y]|&\leq& c(d)\rho^3\left[\theta+
\sigma\rho\Theta^t_T\right](2T+1)^{-d/2},\\
\Theta^t_T&=&\sigma^{-1}\max\limits_{\tau:|\tau|\leq 2T}
|\Delta_1^{\tau_1-t_1}...\Delta_d^{\tau_d-t_d}e|_{2T,\infty}^*.\\
\end{array}
\end{equation}
\end{thm}
The proof of Theorem  \ref{the:Opredict} is identical to that of Theorem \ref{the:Oadapt}.
\section{Families of nice signals}
When applying Algorithms A, B and Theorems \ref{the:Oadapt},
\ref{the:Opredict}, the crucial question is how to recognize
niceness. 
We are
about to give a partial answer to this question.
\subsection{Calculus of nice signals}
\label{calculus} Our current goal is to understand how wide are
the families of nice signals, and our
plan is as follows: (a) we list a number of operations which
preserve the property in question,
 and
 (b) we present a list of examples
of signals possessing the property. Applying to ``raw materials''
from (b) operations from (a), one can produce a wide variety of
nice signals.
Here is a sample of operations preserving niceness of signals.
\\~\\
{\sl I. ``Scale'' of nice signals.} We
start with the following evident observation: $
\rho'\geq\rho,\theta'\geq\theta,L'\leq L\Rightarrow {\cal
F}^t_L(\theta,\rho)\subset{\cal F}^t_{L'}(\theta',\rho') $ and $
\rho'\geq\rho,\theta'\geq\theta,\kappa'\leq\kappa, T_0^\prime\geq
T_0, L'\leq L\Rightarrow {\cal
P}^t_{\kappa,T_0,L}(\theta,\rho)\subset{\cal
P}^t_{\kappa',T_0^\prime,L'}(\theta',\rho'). \\
$ {\sl II. Taking
linear combinations.}
 Our next
observation is that a linear combination of ``good'' signals is
again good, with properly updated parameters:
\begin{prop}\label{propnew} (i) Let $(s^j_\tau)\in{\cal F}^t_L(\theta_j,\rho_j)$, and let $\lambda_j\in{\Bbb C}$ be random
variables independent of $(s^j)$ and such that
$E\{|\lambda_j|^2\}<\infty$, $j=1,...,m$. Then
\begin{equation}\label{eqnewnew1}
\begin{array}{l}
(s_\tau\equiv \sum\limits_{j=1}^m\lambda_j s^j_\tau)\in {\cal
F}^t_{L^+}(\theta^+,\rho^+),\\
\theta^+=(2m-1)^{d/2}2^{m-1}\rho_1...\rho_m\sum\limits_{j=1}^m{\theta_j[E\{|\lambda_j|^2\}]^{1/2}\over
\rho_j},\\
 \rho^+=(2m-1)^{d/2}2^{m}\rho_1...\rho_m,
L^+=\lfloor L/2\rfloor.\\
\end{array}
\end{equation}
In the case of $m=1$, one can  set $\rho^+=\rho_1$,
$\theta^+=|\lambda_1|\theta_1$, $L^+=L$. The filters certifying
the well-filterability of $(s_\tau)$ can be chosen to be
independent of the coefficients $\lambda_j$.
\\
(ii) Let $(s^j_\tau)\in{\cal
P}^t_{\kappa_j,T_0^j,L}(\theta_j,\rho_j)$, $j=1,...,m$, and let
$\lambda_j\in{\Bbb C}$ be random variable independent of $(s^j)$
and such that $E\{|\lambda_j|^2\}<\infty$, $j=1,...,m$. Then
\begin{equation}\label{eqnewnew1pred}
\begin{array}{l}
(s_\tau\equiv \sum\limits_{j=1}^m\lambda_j s^j_\tau)\in {\cal
P}^t_{\kappa^+,L^+}(\theta^+,\rho^+),\\
\theta^+=(2m-1)^{d/2}2^{m-1}\rho_1...\rho_m\sum\limits_{j=1}^m{\theta_j[E\{|\lambda_j|^2\}]^{1/2}\over
\rho_j},\\ \rho^+=(2m-1)^{d/2}2^{m}\rho_1...\rho_m,
\kappa^+=\min\limits_{1\leq j\leq m} \kappa_j,\,
T_0^+=m\max\limits_{1\leq j\leq m}T_0^j,\\
L^+=\lfloor L/2\rfloor.\\
\end{array}
\end{equation}
In the case of $m=1$, one can  set $\rho^+=\rho_1$,
$\theta^+=|\lambda_1|\theta_1$, $\kappa^+=\kappa$, $T_0^+=T_0$,
$L^+=L$. The filters certifying the well-predictability of
$(s_\tau)$ can be chosen to be independent of the coefficients
$\lambda_j$.
\end{prop}
{\sl III. Modulation and conjugation.}  Next we notice
that the families of nice signals are closed w.r.t.
``modulation'' and conjugation:
\begin{prop}\label{propnew1} (i) Let $(s_\tau)\in{\cal F}^t_L(\theta,\rho)$, and let $\omega\in{\Bbb R}^d$,
$\phi\in{\Bbb R}$ be deterministic. Then the signal
$
(\widehat{s}_\tau=\exp\{i[\omega^T\tau+\phi]\}s_\tau)_{\tau\in\bZd}
$
 belongs to ${\cal F}^t_L(\theta,\rho)$ along with $(s_\tau)$, and the signal $(\bar{s}_\tau=\overline{s_\tau})_{\tau}$
($\overline{a}$ is the complex conjugate of $a\in{\Bbb C}$)
belongs to ${\cal F}^t_L(\theta,\rho)$.\\
(ii)
 Let $(s_\tau)\in{\cal P}^t_{\kappa,T_0,L}(\theta,\rho)$, and
let $\omega\in{\Bbb R}^d$, $\phi\in{\Bbb R}$ be deterministic.
Then the signal $
(\widehat{s}_\tau=\exp\{i[\omega^T\tau+\phi]\}s_\tau)_{\tau\in\bZd}
$ also belongs to ${\cal P}^t_{\kappa,T_0,L}(\theta,\rho)$, and
the signal $(\bar{s}_\tau=\overline{s_\tau})_{\tau}$ belongs to
${\cal P}^t_{\kappa,T_0,L}(\theta,\rho)$.
\end{prop}
{\sl IV. Lifting.} We are about to show that a nice
signal in a dimension $d\leq d^+$ can be viewed as a nice signal,
with properly updated parameters, in a dimension $d^+>d$:
\begin{prop}\label{newprop2} (i) Let $1\leq d\leq d^+$, and let $(s_\tau)_{\tau\in\bZd}$
be a signal which is well-filtered, with parameters
$\theta,\rho,L$, at a point $t\in\bZ^{d}$. Then the signal
$(s^+_{\tau_1,...,\tau_{d^+}}=s_{\tau_1,...,\tau_d})$ is
well-filtered, with the parameters $
\theta^+=(2L+1)^{(d^+-d)/2}\theta$, $\rho^+=\rho$, $L^+=L$
 at every point
$t^+\in\bZ^{d^+}$ such that $(t^+_1,...,t^+_d)=t$.
\\
(ii)  Let $1\leq d\leq d^+$, and let $(s_\tau)_{\tau\in\bZd}$ be a
signal which is well-predictable, with parameters
$\theta,\rho,\kappa,T_0,L$, at a point $t\in\bZd$. Then the signal
$(s^+_{\tau_1,...,\tau_{d^+}}=s_{\tau_1,...,\tau_d})$ is
well-predictable, with the parameters
$
\theta^+=(2L+1)^{(d^+-d)/2}\theta,\,\,\rho^+=(2\kappa+1)^{(d_+-d)/2}\rho,\,\,\kappa^+=\kappa,\,\,T_0^+=T_0,\,\,L^+=L,
$
 at every point
$t^+\in\bZ^{d^+}$ such that $(t^+_1,...,t^+_d)=t$.
\end{prop}
{\sl V. ``Tensor product''.} Let $d=d'+d''$ with positive integers
$d'$, $d''$, so that $\bZd=\bZ^{d'}\times\bZ^{d''}$. Given random
fields
 $(s^\prime_{\tau'}(\xi))_{{\tau'\in\bZ^{d'}\atop \xi}}$,
$(s^{\prime\prime}_{\tau''}(\xi))_{{\tau''\in\bZ^{d''}\atop
\xi}}$, we define their tensor product as the field
$
(s_\tau(\xi)=s^\prime_{\tau'}(\xi)s^{\prime\prime}_{\tau''}(\xi))_{{\tau=(\tau',\tau'')\in\bZd\atop
\xi}}.
$
\begin{prop} \label{newprop11} (i) Let $(s^\prime_{\tau'}(\xi))_{{\tau'\in\bZ^{d'}\atop \xi}}
\in {\cal F}^{t'}_{L}(0,\rho')$, $
(s^{\prime\prime}_{\tau''}(\xi))_{{\tau''\in\bZ^{d''}\atop
\xi}}\in {\cal F}^{t''}_{L}(0,\rho'')$. Then $ (s_\tau)\in {\cal
F}^{(t',t'')}_L(0,\rho'\rho''). $
\\
(ii) Let $(s^\prime_{\tau'}(\xi))_{{\tau'\in\bZ^{d'}\atop \xi}}\in
{\cal P}^{t'}_{\kappa,T_0,L}(0,\rho')$,
$(s^{\prime\prime}_{\tau''}(\xi))_{{\tau''\in\bZ^{d''}\atop
\xi}}\in {\cal P}^{t''}_{\kappa,T_0,L}(0,\rho'')$. Then
$(s_\tau)\in {\cal P}^{(t',t'')}_{\kappa,T_0,L}(0,\rho'\rho''). $
\end{prop}
\subsection{Examples of nice signals}
{\sl I. Exponential and algebraic polynomials.} Let us define an
{\sl exponential polynomial} $(s_\tau)$ on $\bZd$ as a {\sl
finite} sum of {\sl exponential monomials}
$c\tau^\alpha\exp\{\omega^T\tau\}\equiv
c\tau_1^{\alpha_1}...\tau_d^{\alpha_d}\exp\{\omega^T\tau\}$ with
{\sl nonnegative} multi-indices $\alpha$ and $\omega\in{\Bbb
C}^d$:
\begin{equation}\label{polyfield}
s_\tau=\sum\limits_{\ell=1}^Mc_\ell\tau^{\alpha(\ell)}\exp\{\omega^T(\ell)\tau\},
\end{equation}
where  $\omega(\ell)$ and $\alpha(\ell)$ are deterministic, and
$c_\ell$ may be random. Given an exponential polynomial $(s_\tau)$
on $\bZd$, we define its {\sl partial sizes} $N_j$, $j=1,...,d$,
as follows: let $m_j$ be the maximum of the degrees
$\alpha_j(\ell)$, $\ell=1,...,M$, of the variable $\tau_j$ in the
monomials of the sum (\ref{polyfield}), and $M_j$ be the number of
{\sl distinct from each other} complex numbers among the ``partial
frequencies'' $\omega_j(\ell)$: $ M_j=\Card{\cal O}_j,\quad {\cal
O}_j=\{\omega_j(\ell): 1\leq \ell \leq M\}. $ The {\sl $j$-th
partial size} $N_j(s)$ of exponential polynomial (\ref{polyfield})
is, by definition, the integer $(m_j+1)M_j$. For example, with all
frequencies equal to 0, an exponential polynomial becomes an
algebraic polynomial, and its $j$-th size is by 1 larger than the
degree of the polynomial w.r.t. $j$-th variable $\tau_j$.
\begin{prop}\label{prophomog} Let $(s_\tau)$ be an exponential polynomial on $\bZd$
of partial sizes $N_1,...,N_d$. Then for all $t\in\bZd$ one has
\begin{equation}\label{expol}
 (s_\tau)\in
{\cal F}^t_\infty(0,\rho_d(N_1,...,N_d)),\, \rho_d(N_1,...,N_d)
=\prod\limits_{j=1}^d[(2N_j-1)^{1/2}2^{3N_j/2}],
\end{equation}
and the filters $q^{(T)}$ certifying this inclusion can be chosen
to be dependent solely on $T$ and on the collection of $d$ sets
${\cal O}_j=\{\omega_j(\ell):1\leq \ell \le M\}$ of partial
frequencies.
\end{prop}
\begin{rem}\label{rempoly} A major shortcoming of {\rm (\ref{expol})}
is a dramatic growth of $\rho_d(N,N,...,N)$ with $N$ and $d$.
In several important cases, better bounds for $\rho$ can be found.
For example, an algebraic polynomial of degree $m$ in every
variable
\begin{equation}\label{polyp}
p_\tau=\sum\limits_{\alpha\geq0,|\alpha|\leq m} c_\alpha
\tau^\alpha
\end{equation}
 belongs to ${\cal F}^t_\infty(0,(16m)^d)$ for every $t$, and the
 filters $q^{(T)}$ certifying this inclusion can be chosen to
 depend solely on $T,d,m$.
 \end{rem}
 {\sl II. Solutions to homogeneous difference equations and
 harmonic functions.} Consider a difference operator ${\cal D}$:
\begin{equation}\label{operator}
({\cal D}f)_\tau=\sum\limits_{\ell=1}^k w_\ell
f_{\tau-\alpha(\ell)};
\end{equation}
here $\alpha(1),...,\alpha(k)\in\bZd$ and $w_1,...,w_k\in{\Bbb
C}$. For a positive integer $N$ and $t\in\bZd$, let $$
\begin{array}{l}B^t_N=\{\tau\in\bZd\mid\,
|\tau-t|\leq N\},\, B^t_N({\cal D})=\{\tau\in
B^t_N\mid\,\tau+\alpha(\ell)\in
B^t_N,\,\ell=1,...,k\},\\
{\cal H}^t_N({\cal D})=\{(s)\in C(\bZd)\mid\,
s_\tau=({\cal D}s)_\tau\quad\forall \tau\in B^t_N({\cal D})\}.\\
\end{array}
$$
For example, with
\begin{equation}\label{eq:average}
({\cal D}f)_\tau={1\over
2d}\sum\limits_{{i=1,...,d\atop\epsilon=\pm1}}
f_{\tau_1,...,\tau_{i-1},\tau_i+\epsilon,\tau_{i+1},...,\tau_d},
\end{equation}
the linear space ${\cal H}^t_N({\cal D})$ is the space of fields
which are ``discrete harmonic'' on $B^t_N$, that is, $
s_\tau={1\over 2d} \sum\limits_{{i=1,...,d\atop\epsilon=\pm1}}
s_{\tau_1,...,\tau_{i-1},\tau_i+\epsilon,\tau_{i+1},...,\tau_d}$
for all $\tau$ with $|\tau-t|\leq N-1. $\\ Let us call a
difference operator ${\cal D}$ {\sl regular}, if it possesses the
following properties:
\\
{\sl R.1} The vectors $\{\alpha(\ell)\}_{1\leq \ell\leq k}$ span
the entire ${\Bbb R}^d$;\\
{\sl R.2} The coefficients $w_\ell=\rho_\ell\exp\{i\phi_\ell\}$
($\rho_\ell\geq0$, $\phi_\ell\in{\Bbb R}$) are nonzero, and
\begin{equation}\label{eq:harmonic1}
\begin{array}{lrclclrcl}
(a)&\sum\limits_{\ell=1}^k\rho_\ell&\leq&1;&\qquad&
(b)&\sum\limits_{\ell=1}^k\rho_\ell\alpha(\ell)&=&0.\\
\end{array}
\end{equation}
For example, the averaging operator (\ref{eq:average}) and its
degrees are regular.
\\
It turns out that the solutions of homogeneous difference
equations with regular difference operators are well-filtered:
\begin{prop}\label{prop:harmonic} Let ${\cal D}$ be a regular
difference operator. Then there exists a constant $c=c({\cal
D})>0$ such that
\begin{equation}
\forall N>0:\quad {\cal H}^t_N({\cal D})\subset {\cal
F}^t_{\lfloor cN\rfloor}(0,c^{-1}).
\end{equation}
\end{prop}
As a nontrivial application example for Proposition
\ref{prop:harmonic}, consider the families of random fields
defined as follows. Let $d\leq 4$, $M$ be a positive integer, and
$R$ be a positive real. Consider the family ${\cal H}^+ (M)$ of
all deterministic continuous functions $f$ on ${\Bbb R}^d$ which
are harmonic in the interior of the box $D^0_{2M}=\{x\in{\Bbb
R}^d:|x_j|\leq 2M, j\leq d\}$:
$
\left(\sum\limits_{j=1}^d{\partial^2\over\partial x_j^2}\right)
f(x)=0,\,x\in\inter D^0_{2M}.
$
Now let ${\cal H}^+(M,R)$ be the family of random functions $f$
such that all realizations of a function belong to ${\cal H}^+(M)$
and, besides this, $E\{\|f\|_{\infty,2M}^2\}\leq R^2$, where
$\|f\|_{\infty,2M}$ is the uniform norm on $D^0_{2M}$. Restricting
functions $f$ from ${\cal H}^+(M,R)$ on $\bZ^d$, we get a family
of random fields ${\cal H}(M,R)$ on $\bZd$.
\begin{prop}\label{propharm} Let $d\leq 4$, $M$ be a positive integer and $R>0$ be a real.
For an appropriately chosen absolute constant $c>0$, for all
deterministic fields $(s_\tau)\in{\cal H}(M,R)$  one has
\begin{equation}\label{harmwellf}
|t|\leq cM,L\leq cM\Rightarrow (s_\tau)\in {\cal
F}^t_L(c^{-1}R,c^{-1}),
\end{equation}
and the filters $q^{(T)}$ certifying the above inclusion can be
chosen depending solely on $d$, $T$.
\end{prop}
\subsection{Basic example of well-predicted signal: quasi-stable exponential polynomial}
Let us define a {\sl quasi-stable} exponential polynomial
$(s_\tau)$ on $\bZd$ as an exponential polynomial
\begin{equation}\label{polyfieldqs}
s_\tau=\sum\limits_{\ell=1}^Mc_\ell\tau^{\alpha(\ell)}\exp\{\omega^T(\ell)\tau\}
\end{equation}
where all partial frequencies $\omega_j(\ell)$ satisfy the
restriction $\Re(\omega_j(\ell))\leq0$.
 For example, an algebraic polynomial (partial frequencies are zero) and  a trigonometric
 polynomial (partial frequencies are imaginary) are quasi-stable.
\begin{prop}\label{prophomogpred} Let $(s_\tau)$ be a quasi-stable exponential polynomial on $\bZd$
of partial sizes $N_1,...,N_d$. Then for every integer
$\kappa\geq0$ and all $t\in\bZd$ one has
\begin{equation}\label{expolpred}
\begin{array}{l}
 (s_\tau)\in
{\cal
P}^t_{\kappa,T_0,\infty}(0,\rho_{\kappa,d}(N_1,...,N_d)),\\
\rho_{\kappa,d}(N_1,...,N_d)
=\prod\limits_{j=1}^d[(2N_j-1)^{1/2}2^{N_j}\left(\max[2,2\kappa+1]\right)^{N_j/2}],\\
T_0=\kappa \max\limits_{1\leq j\leq d} N_j\\
\end{array}
\end{equation}
and the filters $q^{(T)}$ certifying this inclusion can be chosen
to be depending solely on $T,\kappa$ and on the collection of $d$
sets ${\cal O}_j=\{\omega_j(\ell):1\leq \ell \le M\}$ of partial
frequencies.
\end{prop}

\section{Appendix}
\subsection{Preliminaries}
{\sl Norm relations.} Let us list several evident relations
between the introduced semi-norms on $ C({\Bbb  Z}^d)$. \\
$\bullet$ {[Parseval equality]}:
\begin{equation}\label{peq1}
(r,s)_T\equiv \sum\limits_{t:|t|\leq T} r_{t}\overline{s_{t}} =
\sum\limits_{\mu\in\Gamma_T^d}
(F_Tr)(\mu)\overline{(F_Ts)(\mu)}\equiv \langle F_Tr,
F_Ts\rangle_T,
\end{equation}
where $\overline{a}$ is the complex conjugate of $a\in{\Bbb C}$;
in particular,
\begin{equation}\label{peq2}
|r|_{T,2} = |r|_{T,2}^*;
\end{equation}
A useful corollary of Parseval's equality combined with the fact
that $|q|_{T,p}^*=|\bar{q}|_{T,p}^*$ is the relation
\begin{equation}\label{eq:corPars1}
\left|\sum\limits_{|t|\leq T}a_{t}b_{t}\right|\leq
|a|_{T,1}^*|b|_{T,\infty}^*.
\end{equation}
$\bullet$ {[Norms of convolutions of filters]}
\begin{equation}\label{peq3}
r,s\in  C({\Bbb  Z}^d) \Rightarrow
|r(z_1,...,z_d)s(z_1,,...,z_d)|_p\le |r|_1|s|_p;
\end{equation}
$\bullet$ {[Relations between $|\cdot|$ and $|\cdot|^*$]}: for
$p,q\in[1,\infty]$ one has
\begin{equation}\label{peq4}
|r|_{T,p}^* \le
(2T+1)^{d[(1/p-1/2)_++(1/2-1/q)_+]}\,|r|_{T,q},\quad
a_+=\max[a,0];
\end{equation}
\begin{equation}\label{peq5}
\ord(r) +\ord(s) \le T \Rightarrow
|r(z_1,...,z_d)s(z_1,...,z_d)|_{T,p}^*\le |r|_1|s|_{T,p}^*.
\end{equation}
{\sl Useful fact.} In the sequel,  we need the following simple
and well-known fact:
\begin{lem}\label{lem:max} Let
$f_j=\xi_j+i\eta_j$, $0\leq j< N$, be a sequence of $N$ standard
Gaussian complex-valued random variables, not necessarily
independent of each other. Then
\begin{equation}\label{eq:rough}
\begin{array}{l}
[E\{\max\limits_{0\leq j<N}|f_j|^2\}]^{1/2}\leq\sqrt{2\ln  N+2};\\
P\{\max\limits_{0\leq j<N}|f_j|>u+\sqrt{2\ln
N}\}\leq\exp\{-u^2/2\} \quad\forall u\geq 0.\\
\end{array}
\end{equation}
\end{lem}
{\sl Proof.} We have $$\begin{array}{l}\psi(r)
\equiv P\{\max\limits_{0\leq
j<N}|f_j|>r\}\leq\min[1,N\exp\{-r^2/2\}]
 \Rightarrow\\
P\{\max\limits_{0\leq j< N}|f_j|>u+\sqrt{2\ln N}\}\leq
N\exp\{-(u+\sqrt{2\ln N})^2/2\}\leq \exp\{-u^2/2\};\\
 E\{\max\limits_{0\leq j<N}|f_j|^2\}=-\int\limits_0^\infty r^2 d\psi(r)=
 2\int\limits_0^\infty r\psi(r)dr
 \leq2\int\limits_0^{\sqrt{2\ln N}}rdr\\
 +2N\int\limits_{\sqrt{2\ln N}}^\infty
 r\exp\{-r^2/2\}dr
 =2\ln N+2. \;\;\;\;\;\;\;\;\;\qed\\
 \end{array}
 $$
\subsection{Proof of Theorem \protect{\ref{the:Oadapt}}}  W.l.o.g., we may assume that $t=0$. We denote by
$q^*$ the filter associated with $(s_\tau)$ via the description of
the inclusion $(s_\tau)\in{\cal S}^0_{3T}(\theta,\rho,T)$.  Let us
set
\begin{equation}\label{widehatrho}
|q^*|_2=\widehat{\rho}(2T+1)^{-d/2};\quad
{\kappa}=\theta(2T+1)^{-d/2}\qquad
\left[\widehat{\rho}\leq\rho\right],
\end{equation}
so that
\begin{equation}\label{eq:12rho}
\bar{s}=q^*(\Delta)s \Rightarrow \max\limits_{\tau:|\tau|\leq 3T}
E\{\left|s_{\tau}-\bar{s}_{\tau}\right|^2\}\leq \kappa^2.
\end{equation}
 Finally, let
\begin{equation}\label{eq:Theta}
\Theta_T=\max\limits_{\tau:|\tau|\leq 2T}
|\Delta_1^{\tau_1}...\Delta_d^{\tau_d}e|_{2T,\infty}^*,
\end{equation}
and let
 $\widehat{\phi}$ be the optimal solution, used in Algorithm
A, of the optimization problem (\ref{eq:optim}).
\\
{\sl 1$^0$.} We start with simple technical lemma:
\begin{lem}\label{lem:r1}
 Let
$r(z_1,...,z_d)=(q^*(z_1,...,z_d))^2$. Then $r\in C_{2T}(\bZd)$
possesses the following properties:
\begin{eqnarray}
\lefteqn{\label{eq:r1} |r|_2\leq|r|_{2T,1}^*\le
2^{d/2}\widehat{\rho}^2 (2T+1)^{-d/2};}\\
 \label{eq:r101}
\lefteqn{|r|_1\leq \widehat{\rho}^2;}
\\
\label{eq:r30}
\lefteqn{\left[E\{|(1-r(\Delta))s|^2_{2T,2}\}\right]^{1/2}
\le{\kappa} (\widehat{\rho}+1)(4T+1)^{d/2};}\\
\label{eq:20new} \lefteqn{
\begin{array}{rcl}
|(1-r(\Delta))y|^*_{2T,\infty}&\leq&|(1-r(\Delta))s|_{2T,2}
+(1+\widehat{\rho}^2)\Theta_T\\
\end{array}}
\\
\label{eq:r20} \lefteqn{\begin{array}{rcl}
\left[E\left\{\left(|(1-r(\Delta))y|^*_{2T,\infty}
\right)^2\right\}\right]^{1/2}&\le& \sigma
(1+\widehat{\rho}^2)\sqrt{4d\ln (4T+1)+2}\\
&&+{\kappa}(\widehat{\rho}+1)(4T+1)^{d/2}.
\end{array}}\\
\end{eqnarray}
\end{lem}
{\sl Proof.} {\sl (\ref{eq:r1}):} We have
$$
\begin{array}{rcl}
|r|_{2T,1}^*&=&\sum\limits_{\mu\in\Gamma_{2T}^d}{|r(\mu)|\over(4T+1)^{d/2}}
=\sum\limits_{\mu\in \Gamma_{2T}^d}{|q^*(\mu)|^2\over(4T+1)^{d/2}}
=(4T+1)^{d/2}\sum\limits_{\mu\in\Gamma_{2T}^d}\left|{q^*(\mu)\over(4T+1)^{d/2}}\right|^2\\
&=& (4T+1)^{d/2}(|q^*|_{2T,2}^*)^2 =(4T+1)^{d/2}|q^*|_{2T,2}^2\leq
2^{d/2}\widehat{\rho}^2(2T+1)^{-d/2}.\\
\end{array}
$$
Since $|r|_2=|r|_{2T,2}=|r|^*_{2T,2}\le |r|^*_{2T,1}$,
(\ref{eq:r1}) follows.\\
{\sl (\ref{eq:r101}):} We clearly have $|r|_1\leq |q^*|_1^2\leq
((2T+1)^{d/2}|q^*|_2)^2=\widehat{\rho}^2$.\\
{\sl (\ref{eq:r30}):} Let  $h=(1-q^*(\Delta))s$, so that by virtue
of $(s_\tau)\in{\cal S}^0(\theta,\rho,T)$ and in view of the
origin of $q^*$ we have
\begin{equation}\label{eq:wdhdlt}
\max\limits_{\tau:|\tau|\leq 3T} E\{|h_{\tau}|^2\}\leq\kappa^2.
\end{equation}
Setting $g=(1-r(\Delta))s$, we have
$$
\begin{array}{ll}
g_{\tau}=((1+q^*(\Delta))(1-q^*(\Delta))s)_{\tau}=
((1+q^*(\Delta))h)_{\tau} =h_{\tau}+(q^*(\Delta)h)_{\tau}\\
\Rightarrow
|g_{\tau}|\leq|h_{\tau}|+|q^*|_2|\Delta_1^{-\tau_1}...\Delta_d^{-\tau_d}h|_{T,2}\\
\Rightarrow \left(E\{|g_{\tau}|^2\}\right)^{1/2}\leq
\left(E\{|h_{\tau}|^2\}\right)^{1/2}+|q^*|_2\left(\sum\limits_{\tau':|\tau'-\tau|\leq
T}
E\{|h_{\tau-\tau'}|^2\}\right)^{1/2};\\
\end{array}
$$
applying (\ref{eq:wdhdlt}) and taking into account that
$|q^*|_2=\widehat{\rho}(2T+1)^{-d/2}$, we  come to
\begin{equation}\label{eq:bound1}
\max\limits_{\tau:|\tau|\leq 3T} E\{|((1-r(\Delta))s)_{\tau}|^2\}
\leq \left[{\kappa}(\widehat{\rho}+1)\right]^2,
\end{equation}
and (\ref{eq:r30}) follows.\\
{\sl (\ref{eq:20new}), (\ref{eq:r20}):} We have
$$
\begin{array}{l}
|(1-r(\Delta))y|_{2T,\infty}^*\leq |(1-r(\Delta))s|_{2T,\infty}^*
+|(1-r(\Delta))e|_{2T,\infty}^*\\
\leq|(1-r(\Delta))s|_{2T,2}^* +|(1-r(\Delta))e|_{2T,\infty}^*
=|(1-r(\Delta))s|_{2T,2}
+|(1-r(\Delta))e|_{2T,\infty}^*\\
\leq|(1-r(\Delta))s|_{2T,2} +|e|_{2T,\infty}^*+
\sum\limits_{\tau:|\tau|\leq2T}|r_{\tau}||\Delta_1^{\tau_1}...\Delta_d^{\tau_d}e|_{2T,\infty}^*\\
\leq|(1-r(\Delta))s|_{2T,2} +(1+|r|_1)\max\limits_{\tau:|\tau|\leq
2T} |\Delta_1^{\tau_1}...\Delta_d^{\tau_d}e|_{2T,\infty}^*.
\end{array}
$$
The resulting inequality combines with (\ref{eq:r101}) to yield
(\ref{eq:20new}). Further, from the resulting inequality and
(\ref{eq:r30}) it follows that
$$
\begin{array}{l}
\left(E\left\{\left(|(1-r(\Delta))y|_{2T,\infty}^*\right)^2\right\}
\right)^{1/2}\\ \leq\kappa(\widehat{\rho}+1)(4T+1)^{d/2}+
(1+|r|_1)\bigg(E\bigg\{\underbrace{\left(\max\limits_{\tau:|\tau|\leq
2T}
|\Delta_1^{\tau_1}...\Delta_d^{\tau_d}e|_{2T,\infty}^*\right)^2}_{\Theta_T^2}\bigg\}\bigg)^{1/2}
\\  \leq\kappa
(\widehat{\rho}+1)(4T+1)^{d/2}+(1+\widehat{\rho}^2)\left(E\{\Theta_T^2\}\right)^{1/2}\\
\end{array}
$$
(we have used (\ref{eq:r101})). To derive (\ref{eq:r20}) from the
resulting inequality, it remains to note that
\begin{equation}\label{eq:maxnoise}\left(E\{\Theta_T^2\}\right)^{1/2}
\leq\sigma\sqrt{4d\ln(4T+1)+2}.
\end{equation}
Indeed, the coordinates of the Fourier transform of
$\Delta_1^{\tau_1}...\Delta_d^{\tau_d}e$ are, up to factor
$\sigma$, standard complex-valued Gaussian random variables, so
that $\sigma^{-2} \Theta_T^2$ is the maximum of squared modulae of
$(4T+1)^{2d}$ of these variables; therefore
$E\{\Theta_T^2\}\leq\sigma^2(4d\ln(4T+1)+2)$ by Lemma
\ref{lem:max}. \qed \\
{\sl 2$^0$.} We now study the properties of the solution
$\widehat{\phi}$ of  problem (\ref{eq:optim}).
\begin{lem}\label{lem:indst} One has
\begin{eqnarray}
\lefteqn{ |\widehat{\phi}|_{2T,2}\leq2^{d/2}\rho^2(2T+1)^{-d/2};}
\label{eq:13'}
\\
\lefteqn{|(1-\widehat{\phi}(\Delta))e|^*_{2T,\infty}\leq(1+2^d\rho^2)\Theta_T;}
\label{eq:15''}
\\
\lefteqn{\left[
E\left\{\left(|(1-\widehat{\phi}(\Delta))e|^*_{2T,\infty}\right)^2
\right\}\right]^{1/2} \leq \sigma (1+2^d\rho^2)\sqrt{4d\ln
(4T+1)+2};} \label{eq:15'}
\\
\label{eq:16''}
\lefteqn{|(1-\widehat{\phi}(\Delta))s|^*_{2T,\infty}\leq
|(1-r(\Delta))s|_{2T,2}+2(1+2^d\rho^2)\Theta_T;}
\\
\lefteqn{
\begin{array}{rl}
 \left[E\left\{
\left(|(1-\widehat{\phi}(\Delta))s|^*_{2T,\infty}\right)^2\right\}
\right]^{1/2}\le &2\sigma(1+2^d\rho^2)\sqrt{4d\ln (4T+1)+2}\\
&+{\kappa}(\widehat{\rho}+1)(4T+1)^{d/2}.\\ \end{array}}
\label{eq:16'}
\end{eqnarray}
\end{lem}
{\sl Proof.}  {\sl (\ref{eq:13'}):} $
|\widehat{\phi}|_{2T,2}=|\widehat{\phi}|_{2T,2}^*\leq|\widehat{\phi}|_{2T,1}^*\leq2^{d/2}\rho^2(2T+1)^{-d/2},
$ (the concluding inequality comes from the fact that
$\widehat{\phi}$ is feasible for (\ref{eq:optim})).
\\
{\sl (\ref{eq:15''}), (\ref{eq:15'}):} We have
$$
\begin{array}{l}
|(1-\widehat{\phi}(\Delta))e|^*_{2T,\infty}\le
(1+|\widehat{\phi}|_{2T,1})\max\limits_{\tau:|\tau|\leq
2T}|\Delta_1^{\tau_1}...\Delta_d^{\tau_d}e|^*_{2T,\infty}\\
\le\left(1+(4T+1)^{d/2}|\widehat{\phi}|_{2T,2}
\right)\max\limits_{\tau:|\tau|\leq
2T}|\Delta_1^{\tau_1}...\Delta_d^{\tau_d}e|^*_{2T,\infty}
\\
\leq(1+2^d\rho^2) \max\limits_{\tau:|\tau|\leq
2T}|\Delta_1^{\tau_1}...\Delta_d^{\tau_d}e|^*_{2T,\infty}\\
\end{array}
$$
(we have used (\ref{eq:13'})). The resulting inequality implies
that
$$
\begin{array}{l}
\left[E\left\{\left(|(1-\widehat{\phi}(\Delta))e|^*_{2T,\infty}\right)^2
\right\}\right]^{1/2} \leq(1+2^d\rho^2)
\left[E\left\{\max\limits_{\tau:|\tau|\leq 2T}
\left(|\Delta_1^{\tau_1}...\Delta_d^{\tau_d}e|^*_{2T,\infty}\right)^2\right\}\right]^{1/2}\\
\leq(1+2^d\rho^2)\sigma\sqrt{4d\ln(4T+1)+2}\\
\end{array}
$$
(we have used (\ref{eq:maxnoise})).
\\
{\sl (\ref{eq:16'}), (\ref{eq:16'}):} Note that the polynomial $r$
defined in Lemma \ref{lem:r1} is a feasible solution of the
optimization problem (\ref{eq:optim}) by the first relation in
(\ref{eq:r1}), so that the optimal value in the problem does not
exceed $J(r,y^0_{4T})$. It follows that
$$
\begin{array}{rll}
&(a)&J(\widehat{\phi},y^0_{4T})\leq J(r,y^0_{4T})\\
\Rightarrow&(b)&|(1-\widehat{\phi}(\Delta))y|_{2T,\infty}^*\leq
|(1-r(\Delta))y|_{2T,\infty}^*\\
\Rightarrow&(c)&|(1-\widehat{\phi}(\Delta))s|_{2T,\infty}^*\leq
|(1-\widehat{\phi}(\Delta))e|_{2T,\infty}^*+
|(1-r(\Delta))y|_{2T,\infty}^*\\ \Rightarrow& (d)&
\left[E\left\{\left(|(1-\widehat{\phi}(\Delta))s|_{2T,\infty}^*\right)^2
\right\}\right]^{1/2} \leq
\left[E\left\{\left(|(1-\widehat{\phi}(\Delta))e|_{2T,\infty}^*\right)^2
\right\}\right]^{1/2}\\
&&+\left[E\left\{\left(|(1-r(\Delta))y|_{2T,\infty}^*\right)^2
\right\}\right]^{1/2}\\
\end{array}
$$
Relation (\ref{eq:16''}) follows from  $(c)$ combined with
(\ref{eq:20new}) and (\ref{eq:15''}) (recall that
$\widehat{\rho}\leq\rho$). Relation (\ref{eq:16'}) follows from
$(d)$ combined with (\ref{eq:15'}) and (\ref{eq:r20}). \qed \\
{\sl 3$^0$.} Our next step is to prove
\begin{lem}\label{lem:rest}
One has
\begin{eqnarray}
\label{eq:rest2''} \lefteqn{
\begin{array}{rcl}
\left|\left((1-r(\Delta))(1-\widehat{\phi}(\Delta))s\right)_0\right|&\leq&
\left|\left((1-r(\Delta))s\right)_{0}\right|\\
&&+2^{d/2}\rho^2(2T+1)^{-d/2} \left|
(1-r(\Delta))s\right|_{2T,2};\\
\end{array}}\\
\label{eq:rest2}
\lefteqn{\left[E\left\{\left|\left((1-r(\Delta))(1-\widehat{\phi}(\Delta))s\right)_0\right|^2\right\}
\right]^{1/2} \le  {\kappa} (\widehat{\rho}+1)(2^d\rho^2+1).}
\end{eqnarray}
\end{lem}
{\sl Proof.} We have
$$
\begin{array}{lr}
\left|\left((1-r(\Delta))(1-\widehat{\phi}(\Delta))s\right)_{0}
\right|
 \leq\left|\left((1-r(\Delta))s\right)_{0}\right|+
\left|\left(\widehat{\phi}(\Delta)(1-r(\Delta))s\right)_{0}
\right|&\\
\leq\left|\left((1-r(\Delta))s\right)_{0}\right|+
|\widehat{\phi}|_{2T,2}\left| (1-r(\Delta))s\right|_{2T,2}&\\
\leq\left|\left((1-r(\Delta))s\right)_{0}\right|+2^{d/2}\rho^2(2T+1)^{-d/2}
\left| (1-r(\Delta))s\right|_{2T,2}&\hbox{[see (\ref{eq:13'})]}\\
\end{array}
$$
as required in (\ref{eq:rest2''}). From the resulting inequality
it follows that
$$
\begin{array}{lr}
\left[E\left\{\left|\left((1-r(\Delta))(1-\widehat{\phi}(\Delta))s\right)_{0}
\right|^2\right\} \right]^{1/2}
\leq\left[E\left\{\left|\left((1-r(\Delta))s\right)_{0}\right|^2\right\}
\right]^{1/2}&\\
+ 2^{d/2}\rho^2(2T+1)^{-d/2}\left[E\left\{\left|
(1-r(\Delta))s\right|_{2T,2}^2\right\} \right]^{1/2}&\\
\leq{\kappa}(\widehat{\rho}+1)+2^{d/2}\rho^2(2T+1)^{-d/2}
\left[E\left\{\left| (1-r(\Delta))s\right|_{2T,2}^2\right\}
\right]^{1/2}& \hbox{[see
(\ref{eq:bound1})]}\\
\leq{\kappa}(\widehat{\rho}+1)+2^{d/2}\rho^2(2T+1)^{-d/2}
{\kappa}(\widehat{\rho}+1)(4T+1)^{d/2}& \hbox{[see (\ref{eq:r30})]}\\
\end{array}
$$
and (\ref{eq:rest2}) follows. \qed
\\
{\sl 4$^0$.} Now we are able to complete the proof of Theorem
\ref{the:Oadapt}. The error of the estimate $\widehat{s}$ at the
point $t=0$ is
\begin{equation}\label{eq:error}
\begin{array}{l}
s_{0}-\widehat{s}_{0}=
s_{0}-(\widehat{\phi}(\Delta)y)_{0}=\left((1-\widehat{\phi}(\Delta))s\right)_{0}-
\left(\widehat{\phi}(\Delta)e\right)_{0}\equiv
\epsilon^{(1)}_0+\epsilon^{(2)}_0,\\
\epsilon^{(1)}_\tau=\left((1-\widehat{\phi}(\Delta))s\right)_{\tau},\quad
\epsilon^{(2)}_{\tau} = \left(\widehat{\phi}(\Delta)e\right)_{\tau}.\\
\end{array}
\end{equation}
Setting $f_\tau=\overline{e_{-\tau}}$, we have
$$
\begin{array}{lr}
|\epsilon^{(2)}_{0}|=\left|\sum\limits_{\tau:|\tau|\leq 2T
}\widehat{\phi}_\tau e_{-\tau}\right| \leq
|\widehat{\phi}|_{2T,1}^*|f|_{2T,\infty}^*&\hbox{[see
(\ref{eq:corPars1})]}\\ \leq
2^{d/2}\rho^2(2T+1)^{-d/2}|f|_{2T,\infty}^*,&\hbox{[since
$\widehat{\phi}$ is feasible for
(\ref{eq:optim})]}\\
\end{array}
$$
whence, by definition of $\Theta_T$,
\begin{equation}\label{eq:epsilon2_0}
|\epsilon^{(2)}_0|\leq2^{d/2}\rho^2(2T+1)^{-d/2}\Theta_T.
\end{equation}
Applying (\ref{eq:maxnoise}), we derive from the latter inequality
that
\begin{equation}\label{eq:error1}
\left[E\left\{|\epsilon^{(2)}_{0}|^2\right\}\right]^{1/2} \leq
2^{d/2}\sigma\rho^2(2T+1)^{-d/2} \sqrt{2d\ln(4T+1)+2}.
\end{equation}
We further have
\begin{eqnarray}
|\epsilon^{(1)}_{0}|&=&\left|\left((1-\widehat{\phi}(\Delta))s\right)_{0}\right|\nn
&\leq&\left|\left(r(\Delta)(1-\widehat{\phi}(\Delta))s\right)_{0}\right|
+
\left|\left((1-r(\Delta))(1-\widehat{\phi}(\Delta))s\right)_{0}\right|\nn
&\underbrace{\leq}_{\small a}&
|r|_{2T,1}^*|(1-\widehat{\phi}(\Delta))s|_{2T,\infty}^*
+\left|\left((1-r(\Delta))(1-\widehat{\phi}(\Delta))s\right)_{0}\right|
\nn &\underbrace{\leq}_{\small b}&
2^{d/2}\rho^2(2T+1)^{-d/2}|(1-\widehat{\phi}(\Delta))s|_{2T,\infty}^*\nn &&
+\left|\left((1-r(\Delta))(1-\widehat{\phi}(\Delta))s\right)_{0}\right|
\label{wenowave}
\end{eqnarray}
(the inequality $a$ is given by (\ref{eq:corPars1}), and  $b$ follows from the feasibility of
$\widehat{\phi}$ for (\ref{eq:optim})), whence
\begin{equation}\label{eq:error2}
\begin{array}{l}
\left[E\left\{|\epsilon^{(1)}_{0}|^2\right\}\right]^{1/2}
\leq2^{d/2}\rho^2(2T+1)^{-d/2}
\left[E\left\{\left(\left|(1-\widehat{\phi}(\Delta))s\right|_{2T,\infty}^*
\right)^2\right\}\right]^{1/2}\\
+
\left[E\left\{\left|\left((1-r(\Delta)(1-\widehat{\phi}(\Delta))s\right)_{0}
\right|^2\right\} \right]^{1/2}\\
\leq2^{d/2}\rho^2(2T+1)^{-d/2}\bigg[2^d\sigma(1+2^d\rho^2)\sqrt{4d\ln
(4T+1)+2}\\ +{\kappa}(\widehat{\rho}+1)(4T+1)^{d/2}\bigg]
+{\kappa} (\widehat{\rho}+1)(2^d\rho^2+1)\\
\end{array}
\end{equation}
(see (\ref{eq:16'}), (\ref{eq:rest2})). Combining
(\ref{eq:error}), (\ref{eq:error1}), (\ref{eq:error2}), we finally
get
\begin{equation}\label{eq:final}
\begin{array}{l}
\left[E\left\{\left|s_{0}-\widehat{s}_{0}\right|^2\right\}\right]^{1/2}
\leq 2^{d/2}\sigma\rho^2(2T+1)^{-d/2}\sqrt{2d\ln(4T+1)+2}\\
+
2^{d/2}\rho^2(2T+1)^{-d/2}\bigg[2^d\sigma(1+2^d\rho^2)\sqrt{4d\ln
(4T+1)+2}\\ +{\kappa}(\widehat{\rho}+1)(4T+1)^{d/2}\bigg]
+{\kappa}
(\widehat{\rho}+1)(2^d\rho^2+1).\\
\end{array}
\end{equation}
Recalling that $\widehat{\rho}\leq\rho$,
$\kappa=\theta(2T+1)^{-d/2}$ and that $\rho\geq1$, (\ref{eq:0th})
follows. \\  Now assume that $(s)$ is deterministic. In this case,
from (\ref{wenowave}) combined with (\ref{eq:16''}) and
(\ref{eq:rest2''}) implies that
\begin{equation}\label{howhowIamtired}
\begin{array}{l}
|\epsilon^{(1)}_{0}|\leq2^{1+d/2}\rho^2(2T+1)^{-d/2}|(1-r(\Delta))s|_{2T,2}\\
+2^{1+d/2}\rho^2(1+2^d\rho^2)(2T+1)^{-d/2}\Theta_T
+\left|\left((1-r(\Delta))s\right)_{0}\right|,
\end{array}
\end{equation}
while from (\ref{eq:r30}), (\ref{eq:bound1}) it follows that
\begin{equation}\label{Iamclose}
\begin{array}{l}
|(1-r(\Delta))s|_{2T,2}\leq{\kappa}
(\widehat{\rho}+1)(4T+1)^{d/2}\leq 2^{d/2}\theta(1+\rho),\\
|((1-r(\Delta))s)_0|\leq\kappa(1+\widehat{\rho})\leq
\theta(1+\rho)(2T+1)^{-d/2}.\\
\end{array}
\end{equation}
 Therefore (\ref{howhowIamtired}) implies that
\begin{equation}\label{IVIV}
|\epsilon^{(1)}_{0}|\leq
3^{3+d}\rho^3\left[\theta+\rho\Theta_T\right](2T+1)^{-d/2}.
\end{equation}
Combining this relation with (\ref{eq:epsilon2_0}) and
(\ref{eq:error}), we arrive at (\ref{ThatIsIt}). \qed
%
%
%
%
%
%
\subsubsection{Proof of Proposition \protect{\ref{propnew}}}
In the proofs to follow, we focus on the
case of well-filtered signals; the reasoning in the case of
well-predicted signals is completely similar.
\par
The
case of $m=1$ is evident. Now let $m\geq2$, let $T^+$ be an
integer, $0\leq T^+\leq L^+$, and let $T=\lfloor
m^{-1}T^+\rfloor$. Since $s^j\in{\cal F}^t_L(\theta,\rho)$ and
clearly $T\leq L$, there exist filters $q^j$ such that
\begin{equation}\label{eqnewnew2}
\begin{array}{l}
(a):\ord(q^j)\leq T;\,\,(b):|q^j|_2\leq\rho_j(2T+1)^{-d/2};\\
(c):|q^j|_1=|q^j|_{T,1}\leq (2T+1)^{d/2}|q^j|_2
\leq\rho_j;\\
(d):
\left[E\left\{|s^j_\tau-(q^j(\Delta)s^j)_\tau|^2\right\}\right]^{1/2}\leq
\theta_j(2T+1)^{-d/2}\quad \forall (\tau:|\tau-t|\leq L).\\
\end{array}
\end{equation}
Now let filter $q$ be defined by $$
1-q(z)=\prod\limits_{j=1}^m(1-q^j(z)),\,\,z=(z_1,...,z_d). $$
Observe that
\begin{equation}\label{eqnewnew3}
\ord(q)\leq mT\leq T^+.
\end{equation}
Note that \begin{equation}\label{eqnewnew4} |q|_2\leq
2^m\rho_1...\rho_m (2T+1)^{-d/2}\leq
(2m-1)^{d/2}2^m\rho_1...\rho_m(2T^++1)^{-d/2}.
\end{equation}
Indeed, we clearly have
$$
\begin{array}{l}
|q(z)|_2=\left|\sum\limits_{\ell=1}^m(-1)^{\ell+1}\sum\limits_{1\leq
j_1<j_2<...<j_\ell\leq
m}q^{j_1}(z)q^{j_2}(z)...q^{j_\ell}(z)\right|_2\\
\leq\sum\limits_{\ell=1}^m\sum\limits_{1\leq
j_1<j_2<...<j_\ell\leq m}|q^{j_1}(z)q^{j_2}(z)...q^{j_\ell}(z)|_2
\underbrace{\leq}_{a} \sum\limits_{\ell=1}^m\sum\limits_{1\leq
j_1<j_2<...<j_\ell\leq
m}{\rho_{j_1}\rho_{j_2}...\rho_{j_\ell}\over(2T+1)^{d/2}}\\
\leq [(1+\rho_1)...(1+\rho_m)-1](2T+1)^{-d/2}
\underbrace{\leq}_{b}2^m\rho_1...\rho_m(2T+1)^{-d/2}\\
\end{array}
$$
($a$ is by (\ref{eqnewnew2}.$b-c$) since
$|u(z)v(z)|_2\leq|u|_1|v|_2,$ $|u(z)v(z)|_1\leq|u|_1|v|_1$], $b$
is  due to $\rho_j\geq1$), as required in (\ref{eqnewnew4}).
Further, by (\ref{eqnewnew2}.$c$), for the filters
$$
Q^j(z)=\left(\prod\limits_{\ell=1}^{j-1}(1-q^\ell(z))\right)
\left(\prod\limits_{\ell=j+1}^{m}(1-q^\ell(z))\right)
$$
one has
\begin{equation}\label{eqnewnew45}
|Q^j|_1\leq
(1+\rho_1)...(1+\rho_{j-1})(1+\rho_{j+1})...(1+\rho_m)\leq
{2^{m-1}\rho_1...\rho_m\over\rho_j}.
\end{equation}
 Now let
$\tau\in\bZd$ be such that $|\tau-t|\leq L^+$.  We have
\begin{eqnarray*}
\lefteqn{\left[E\left\{\left|\left(1-q(\Delta))s\right)_\tau\right|^2\right\}\right]^{1/2}
=\left[E\left\{\left|\sum\limits_{j=1}^m\lambda_j\left(1-q(\Delta))s^j\right)_\tau\right|^2\right\}\right]^{1/2}
}\\
&\leq&
\sum\limits_{j=1}^m\left[E\left\{\left|\lambda_j\left(1-q(\Delta))s^j\right)_\tau\right|^2\right\}\right]^{1/2}
\\&\underbrace{=}_{a}& \sum\limits_{j=1}^m[E\{|\lambda_j|^2\}]^{1/2}
\left[E\left\{\left|\left(1-q(\Delta))s^j\right)_\tau\right|^2\right\}\right]^{1/2}\\
  &\leq& \sum\limits_{j=1}^m[E\{|\lambda_j|^2\}]^{1/2}
\left[E\left\{\left|\left(Q^j(\Delta)(1-q^j(\Delta))s^j\right)_\tau\right|^2\right\}\right]^{1/2}\\
 &\underbrace{\leq}_{b}& \sum\limits_{j=1}^m[E\{|\lambda_j|^2\}]^{1/2}|Q^j|_1\max\limits_{\tau':|\tau'-\tau|\leq(m-1)T}
\left[E\left\{\left|\left((1-q^j(\Delta))s^j\right)_{\tau'}\right|^2\right\}\right]^{1/2}\\
 &\leq&
 2^{m-1}\rho_1...\rho_m(2T+1)^{-d/2}\sum\limits_{j=1}^m{\theta_j[E\{|\lambda_j|^2\}]^{1/2}\over\rho_j}\\
 &\leq&
 \left[(2m-1)^{d/2}2^{m-1}\rho_1...\rho_m\sum\limits_{j=1}^m{\theta_j[E\{|\lambda_j|^2\}]^{1/2}\over\rho_j}\right]
(2T^++1)^{-d/2}
\end{eqnarray*}
where $a$ is due to independence of $\lambda_j$ and $(s^j)$ and
$b$ follows from (\ref{eqnewnew45}), (\ref{eqnewnew2}.$d$), and since
$$|\tau'-\tau|\leq (m-1)T,|\tau-t|\leq L^+\Rightarrow|\tau'-t|\leq
L^++T^+\leq L.$$ Combining the resulting inequality,
(\ref{eqnewnew3}), (\ref{eqnewnew4})  and taking into account that
$T^+\in\{0,1,...,L^+\}$ is arbitrary, we conclude that $s\in {\cal
F}^t_{L^+}(\theta^+,\rho^+)$. Note that by construction, the
filters certifying the latter inclusion are independent of
$\lambda_j$. \qed \subsubsection{Proof of Proposition
\protect{\ref{propnew1}}} (i): Let $T\leq L$, and let $q$ be such
that
\begin{equation}\label{eqnewnew5}
\begin{array}{l}
\ord(q)\leq T,\, |q|_2\leq{\rho\over(2T+1)^{d/2}},\\\
 \max\limits_{\tau:|\tau-t|\leq L}
\left[E\left\{\left|\left((1-q(\Delta))s\right)_\tau\right|^2\right\}\right]^{1/2}\leq
{\theta\over(2T+1)^{d/2}}.\\
\end{array}
\end{equation}
Let us set
$
\widehat{q}_\tau=\exp\{i\omega^T\tau\}q_\tau,\,\,\tau\in\bZd.
$
Then $\ord(\widehat{q})\leq T$, $|\widehat{q}|_2=|q|_2$ and
$$
\begin{array}{rcl}
\left((1-\widehat{q}(\Delta))\widehat{s}\right)_\tau&=&
\exp\{i[\omega^T\tau+\phi]\}s_\tau\\
&&-\sum\limits_{\tau'}
(\exp\{i\omega^T\tau'\}q_{\tau'})(\exp\{i[\omega^T(\tau-\tau')+\phi]\}s_{\tau-\tau'})\\
&=&
\exp\{i[\omega^T\tau+\phi]\}\left((1-q(\Delta))s\right)_\tau,\\
\end{array}
$$
so that (\ref{eqnewnew5}) remains valid when $q,(s)$ are replaced
with $\widehat{q},(\widehat{s})$. Thus, $(\widehat{s})\in{\cal
F}^t_L(\theta,\rho)$. (i) is proved;  (ii) is evident. \qed
\subsubsection{Proof of Proposition \protect{\ref{newprop2}}} Let
$T\leq L$, and let $q=(q_\tau)_{\tau\in\bZd}$ be such that
$\ord(q)\leq T$, $|q|_2\leq \rho(2T+1)^{-d/2}$,
$$
\left[E\left\{\left|\left((1-q(\Delta))s\right)_\tau\right|^2\right\}\right]^{1/2}
\leq \theta (2T+1)^{-d/2}\quad \forall (\tau\in\bZd:|\tau-t|\leq
L).
$$
Setting
$
q^+_{\tau_1,...,\tau_{d^+}}=(2T+1)^{-(d^+-d)}q_{\tau_1,...,\tau_d},
$
we clearly have $\ord(q^+)\leq T$, $|q^+|_2\leq
\rho(2T+1)^{-d^+/2}$ and
$$
\left[E\left\{\left|\left((1-q^+(\Delta))s^+\right)_\tau\right|^2\right\}\right]^{1/2}
\leq \theta (2T+1)^{-d/2}\quad \forall
(\tau\in\bZ^{d^+}:|\tau-t^+|\leq L).
$$
It remains to note that $\theta(2T+1)^{-d/2}\leq
\theta^+(2T+1)^{-d^+/2}$ for $0\leq T\leq L$. \qed
\subsubsection{Proof of Proposition \protect{\ref{newprop11}}} Let
$T\leq L$, and let $q'\in C_T(\bZ^{d'})$, $q''\in C_T(\bZ^{d''})$
be such that
\begin{equation}\label{asnem112}
\begin{array}{l}
(a):|q'|_2\leq\rho'(2T+1)^{-d'/2},\,|q''|_2\leq\rho''(2T+1)^{-d''/2},\\
s^\prime_{\tau'}(\xi)=\sum\limits_{\nu'}s^\prime_{\tau'-\nu'}q^\prime_{\nu'},\,|\tau'-t'|\leq
L,\\
s^{\prime\prime}_{\tau''}(\xi)=\sum\limits_{\nu''}s^{\prime\prime}_{\tau''-\nu''}q^{\prime\prime}_{\nu''},\,|\tau''-t''|\leq L\\
\end{array}.
\end{equation}
Let
$
q(z_1,...,z_d)=q'(z_1,...,z_{d'})q''(z_{d'+1},...,z_{d}),
$
so that
\begin{equation}\label{asnem114}
q\in C_T(\bZd),\quad |q|_2=|q'|_2|q''|_2\leq
\rho'\rho''(2T+1)^{-d/2}
\end{equation}
(see (\ref{asnem112}.$a$)). Now let $\tau=(\tau',\tau'')$ be such
that $|\tau-(t',t'')|\leq L$. We have
$$
\begin{array}{l}
(q(\Delta)s(\xi))_\tau=\sum\limits_{(\nu',\nu'')\in\bZ^{d'}\times\bZ^{d''}}
s^\prime_{\tau'-\nu'}(\xi)s^{\prime\prime}_{\tau''-\nu''}(\xi)q^\prime_{\nu'}q^{\prime\prime}_{\nu''}
=\sum\limits_{\nu'\in
\bZ^{d'}}s^\prime_{\tau'-\nu'}q^\prime_{\nu'}s^{\prime\prime}_{\tau''}(\xi)\\
=s^\prime_{\tau'}(\xi)s^{\prime\prime}_{\tau''}(\xi)
=s_{(\tau',\tau'')},\\
\end{array}
$$
which combines with (\ref{asnem114}) to yield that
$(s_\tau)\in{\cal F}^{(t',t'')}_L(0,\rho'\rho'')$. \qed
\subsubsection{Proof of Proposition \protect{\ref{prophomog}}} We
start with the following two evident facts:
\begin{lem}\label{lem:ev1} Let $(s^j)\in C(\bZd)$ be deterministic
fields belonging to ${\cal F}^t_L(\theta,\rho)$, $j=1,2,...$ such
that $s^j_\tau\to s_\tau$, $j\to\infty$, for every $\tau\in\bZd$.
Then $(s)\in {\cal F}^t_L(\theta,\rho)$.
\end{lem}
Indeed, for every $T$, $0\leq T\leq L$, the filters $q^{j,T}\in
C_T(\bZd)$ which certify the inclusions $(s^j)\in{\cal
F}^t_L(\theta,\rho)$ satisfy $|q^{j,T}|_2\leq \rho(2T+1)^{-d/2}$
and therefore have a limiting point $q^T\in C_T(\bZd)$ with
$|q^j|_2\leq \rho(2T+1)^{-d/2}$. The filters $\{q^T\}_{0\leq T\leq
L}$ clearly certify the inclusion $(s)\in {\cal
F}^t_L(\theta,\rho)$. \qed
\begin{lem}\label{lem:ev2} For every $t\in\bZ$, the univariate exponential field
$(s_\tau=\exp\{\omega \tau\})$, $\omega\in{\Bbb C}$, belongs to
${\cal F}^t_\infty(0,\sqrt{2})$.\end{lem} Indeed, assuming
$\Re(\omega)\geq 0$ and given $T\geq0$, let us set $q(z)={1\over
T+1}[1+\exp\{-\omega\}z^{-1}+\exp\{-2\omega\}z^{-2}+...+\exp\{-T\omega\}z^{-T}]$.
Then $q\in C_T(\bZ)$, $|q|_2=(T+1)^{-1/2}\leq
2^{1/2}(2T+1)^{-1/2}$, while clearly $q(\Delta)s\equiv s$. In the
case of $\Re(\omega)<0$, the same reasoning holds true for
$q(z)={1\over
T+1}[1+\exp\{\omega\}z+\exp\{2\omega\}z^2+...+\exp\{T\omega\}z^T]$.
\qed\\ To complete the proof, we need the following fact:
\begin{lem}\label{lem:ev3} Let $(s_\tau)$ be a ``simple'' exponential polynomial -- a deterministic
exponential polynomial of the form
$
(s_\tau)=\sum\limits_{\ell=1}^M c_\ell \exp\{\omega^T(\ell)\tau\}.
$
Then
\begin{equation}\label{asnem120}
\forall t\in \bZd: (s_\tau)\in{\cal
F}^t_\infty(0,\rho_d(N_1,...,N_d)),
\end{equation}
where $\rho_d(\cdot,...,\cdot)$ is given by (\ref{expol}) and
$N_1,...,N_d$ are the partial sizes of the polynomial. Besides
this, the filters $q^{(T)}$ certifying the above inclusion can be
chosen to depend solely on $T$ and on the collection of the $d$
sets ${\cal O}_j=\{\omega_j(\ell): \ell=1,...,M\}$.
\end{lem}
{\sl Lemma \protect{\ref{lem:ev3}} $\Rightarrow$ Proposition
\ref{prophomog}:} Assume first that the coefficients $c_\ell$ in
(\ref{polyfield}) are deterministic. Since every one of the
univariate functions $f(t)=t^k$, $0\leq k\leq m$, is, uniformly on
compact sets, the limit, as $\epsilon\to+0$, of appropriate linear
combinations of the $m+1$ exponents $\exp\{-k\epsilon\}$, the
exponential polynomial (\ref{polyfield}) is the pointwise, on
$\bZd$, limit, as $i\to\infty$, of simple exponential polynomials
$(s_\tau^i)$ with extended sets of ``frequencies''
$\{\omega_j(\ell)\}_{j,\ell}$: in the approximating polynomials,
every one of these frequencies is replaced by $(m_j+1)$
frequencies $\omega_j(\ell)-k\epsilon_i$, $0\leq k\leq m_j$. Note
that by the definition of partial sizes of exponential
polynomials, the approximating polynomials have exactly the same
partial sizes as the original polynomial $(s_\tau)$. Combining
Lemmas \ref{lem:ev3} and (\ref{lem:ev1}), we immediately conclude
that the exponential polynomial (\ref{polyfield}) belongs to
${\cal F}^t_\infty(0,\rho_d(N_1,...,N_d))$. Since the filters
$q^{(T),i}$ certifying well-filterability of the approximating
polynomials $(s_\tau^i)$ can be chosen to depend solely on $T$ and
the sets of partial frequencies of these approximating
polynomials, from the proof of Lemma \ref{lem:ev1} it follows that
the filters $q^{(T)}$ certifying the inclusion $(s_\tau)\in {\cal
F}^t_\infty(0,\rho_d(N_1,...,N_d))$ can be chosen to depend solely
on $T$ and the sets of partial frequencies of $(s_\tau)$, as
required in Proposition \ref{prophomog}. We have proved
Proposition \ref{prophomog} for the case of a deterministic
exponential polynomial; since the filters certifying
well-filterability of such a polynomial are independent of the
coefficients $c_\ell$, the result is valid for random polynomials
as well. \qed\\
{\sl Proof of Lemma \ref{lem:ev3}.} Proof is by induction in $d$.
\\
\underline{\sl Base $d=1$} is readily given by Lemma \ref{lem:ev2}
combined with Proposition \ref{propnew}.\\ \underline{\sl Step
$1\leq d\Rightarrow d+1$:} Let $ s_\tau=\sum\limits_\ell c_\ell
\exp\{\omega^T(\ell)\tau\} $ be a simple exponential polynomial on
$\bZ^{d+1}$ with partial sizes $N_j$ and the sets of partial
frequencies ${\cal O}_j$, $j=1,...,N$. Let $T\geq0$, and let
$t\in\bZ^{d+1}$. By the inductive hypothesis, there exist filters
$g^{(T)}\in C_T(\bZd)$, $h^{(T)}\in C_T(\bZ)$ (depending solely on
$T$ and on ${\cal O}_1$,...,${\cal O}_{d+1}$) such that
\begin{equation}\label{asnem121}
\begin{array}{l}
(a):|g^{(T)}|_2\leq \rho_d(N_1,...,N_d)(2T+1)^{-d/2},\\
(a'):|h^{(T)}|_2\leq \rho_1(N_{d+1})(2T+1)^{-1/2},\\
(b):r_\tau =\sum\limits_{\nu\in\bZd} r_{\tau-\nu}g^{(T)}_\nu\,\, \forall \tau\in\bZd\,\,\forall (r_\tau)\in {\cal E}({\cal O}_1,...,{\cal O}_d),\\
(b'):p_\tau =\sum\limits_{\nu\in\bZ} p_{\tau-\nu}h^{(T)}_\nu\,\,\forall \tau\in\bZ\,\,\forall (p_\tau)\in {\cal E}({\cal O}_{d+1}),\\
\end{array}
\end{equation}
where ${\cal E}({\cal O}^1,...,{\cal O}^m)$ is the space of all
simple exponential polynomials on $\bZ^{m}$ with the sets of
partial frequencies ${\cal O}^1,...,{\cal O}^m$. Setting
$
q^{(T)}_\tau=g^{(T)}_{\tau_1,...,\tau_d}h^{(T)}_{\tau_{d+1}},\quad
\tau\in\bZ^{d+1},
$
we clearly have
\begin{equation}\label{asnem123}
\begin{array}{l}
q^{(T)}\in C_T(\bZ^{d+1}),\quad |q^{(T)}|_2=|g^{(T)}|_2|h^{(T)}|_2
\leq
\rho_d(N_1,...,N_d)\rho_1(N_{d+1})\\
=\rho_{d+1}(N_1,...,N_{d+1})\\
\end{array}
\end{equation}
(see (\ref{asnem121}.$a,a'$). Further, for every $(s_\tau)\in{\cal
E}({\cal O}_1,...,{\cal O}_{d+1})$ we have, setting
$\tau=(\tau',\tau'')$ with $\tau'\in\bZd$, $\tau''\in\bZ$:
$$
\begin{array}{l}
\sum\limits_{\nu\in\bZ^{d+1}}q^{(T)}_\nu s_{\tau-\nu}=
\sum\limits_{\nu'\in\bZd}g^{(T)}_{\nu'}\left(
\sum\limits_{\nu''\in
\bZ}h^{(T)}_{\nu''}s_{\tau'-\nu',\tau''-\nu''}\right)
\underbrace{=}_{a}\sum\limits_{\nu'\in\bZd}g^{(T)}_{\nu'}s_{\tau'-\nu',\tau''}\\
\underbrace{=}_{b}s_{\tau',\tau''}\\
 \end{array}
 $$
 ($a$ is by (\ref{asnem121}.$b'$) since
 $(s_{\tau'-\nu',\mu})_{\mu\in\bZ}\in{\cal E}({\cal O}_{d+1})$,
 $b$ is by (\ref{asnem121}.$b$) since
 $(s_{\mu,\tau''})_{\mu\in\bZd}\in{\cal E}({\cal O}_1,...,{\cal
 O}_{d})$),
 which combines with (\ref{asnem123}) to imply that
 $$(s_\tau)\in{\cal
 S}^t_\infty(0,\rho_{d+1}(N_1,...,N_{d+1}),T).$$
 Thus, the filters $q^{(T)}$ (which depend
 solely on $T$ and ${\cal O}_1$,...,${\cal O}_{d+1}$) certify the
 inclusion $(s_\tau)\in{\cal
 L}^t_\infty(0,\rho_{d+1}(N_1,...,N_{d+1})$. The inductive step is
 completed.
 \qed
\subsubsection{Proof of statement in Remark
\protect{\ref{rempoly}}} It suffices to prove that for every
nonnegative integer $T$ and every $m,d$ there exists a filter
$q^{(T)}$, $\ord(q^{(T)})\leq T$, depending solely on $T,m,d$,
such that
\begin{equation}\label{eqeqeq78}
\begin{array}{lrcl}
(a)&q^{(T)}(\Delta)p&=&p \hbox{\ for every polynomial
(\ref{polyp})},\\
(b)&|q^{(T)}|_2&\leq&  \left({16m\over\sqrt{2T+1}}\right)^d\equiv \Theta^d.\\
\end{array}
\end{equation}
This well-known fact can be proved by induction in $d$ completely
similar to the one used to prove Lemma \ref{lem:ev3}; the only
difference is in the Base, which now should be replaced with the
following statement:
\begin{lem}\label{lem:ev4} Let
$ p(\tau)=\sum\limits_{\ell=0}^m p_\ell \tau^\ell $ be a
deterministic univariate algebraic polynomial of degree $m$. Then
for every $T\geq0$ there exists a filter $q\in C_T(\bZ)$,
depending solely on $T,m$, with $|q|_2\leq 16m(2T+1)^{-1/2}$ such
that $p(t)=\sum\limits_\nu q^{(T)}_\nu p(t-\nu)$ for all
$t\in\bZ.$
\end{lem}
{\sl Proof.} By evident reasons, it suffices to prove that for a
given $T\geq0$ there exists a collection of weights $q_t$, $-T\leq
t\leq T$, such that
$$
\sum\limits_{t=-T}^T q_t=1,\, \sum\limits_{t=-T}^Tq_t
t^i=0,\,i=1,...,m,\, \sum\limits_{t=-T}^T q_t^2\leq \Theta^2\equiv
{256 m^2\over 2T+1}.
$$
By the standard separation arguments, this is the same as to prove
that for every real algebraic polynomial $r(t)$ of degree $\leq m$
such that $r(0)=1$ one has
$
\sum\limits_{t=-T}^T r^2(t)\geq {2T+1\over 256m^2},
$
or, which is the same, that for the real trigonometric polynomial
$\rho(\phi)=r(T\sin(\phi))$ one has
\begin{equation}\label{eqeqeq79}
\sum\limits_{t=-T}^T \rho^2(\phi_t) \geq {2T+1\over 256 m^2},\quad
\phi_t = \asin(t/T).
\end{equation}
Note that the degree of the trigonometric polynomial $\rho(\cdot)$
is $\leq m$ and that $\rho(0)=1$. Besides this,
$\rho(\phi)=\rho(\pi-\phi)$; due to the latter fact,
$$
M\equiv \max\limits_\phi |\rho(\phi)|
=\max\limits_{|\phi|\leq{\pi\over 2}}|\rho(\phi)|\geq |\phi(0)|=1.
$$
By Bernstein's Theorem on trigonometric polynomials, we have
$|\rho'(\phi)|\leq mM$. Now let $\phi_*\in[-\pi/2,\pi/2]$ be a
point such that $|\rho(\phi_*)|=M$, let $\widehat{\Delta}$ be the
segment of the length ${1\over m}$ centered at $\phi_*$, and
$\Delta$ be the part of this segment in $[-\pi/2,\pi/2]$. Note
that the length of $\Delta$ is at least ${1\over 2m}$ and that for
$\phi\in\Delta$ one has $|\rho(\phi)|\geq |\rho(\phi_*)|-{1\over
2m}(mM)\geq M/2$. Let $n$ be the minimum number of points $\phi_t$
belonging to a segment $\delta\subset[-\pi/2,\pi/2]$ of the length
$1/(2m)$, the minimum being taken over all positions of $\delta$
in $[-\pi/2,\pi/2]$. It is immediately seen that
$
n\geq (1-\sin(\pi/2-1/(2m)))T-2 \geq {T\over 16m^2}-2,
$
whence
$$
\sum\limits_{t=-T}^T\rho^2(\phi_t)\geq
\sum\limits_{t:\phi_t\in\Delta} \rho^2(\phi_t)\geq {M^2\over
4}n\geq {M^2\over 4}\left[{T\over 16 m^2}-2\right]\geq {1\over
4}\left[{T\over 16 m^2}-2\right].
$$
When $T\geq 64m^2$, the latter quantity is $\geq {2T+1\over
256m^2}$, and in any case
$\sum\limits_{t=-T}^T\rho^2(\phi_t)\geq\rho^2(\phi_0)=1$. Thus, we
always have $ \sum\limits_{t=-T}^T\rho^2(\phi_t)\geq {2T+1\over
256 m^2}, $ as required in (\ref{eqeqeq79}). \qed
\subsubsection{Proof of Proposition \protect{\ref{prop:harmonic}}}
In the proof to follow, $c_i$ stand for positive constants
depending solely on ${\cal D}$.
\\
{\sl 1$^0$.} We start with the following evident observation:
\begin{lem}\label{lem:harm1} There exists $c_1$ such that for
every polynomial $p(t)$ of one variable satisfying the relation
$p(1)=1$ one has
\begin{equation}\label{eq:harm2}
\begin{array}{l}
M\leq c_1N,\, \deg(p)\leq c_1N,\\ (s)\in{\cal H}^t_N({\cal
D})\Rightarrow s_\tau=(p({\cal D})s)_\tau \,\, \forall
(\tau:|\tau-t|\leq M).
\\
\end{array}
\end{equation}
\end{lem}
{\sl 2$^0$.} Let us fix a positive integer $N$, and let
\begin{equation}\label{deltaandOmega}
\begin{array}{l}
\delta(\omega)=\sum\limits_{\ell=1}^k w_\ell
\exp\{i\omega^T\alpha(\ell)\}: [-\pi,\pi]^d\to{\Bbb C},\\
\Omega_N^d=\left\{\omega\in{\Bbb R}^d\mid\, \omega_j\in
\left\{{q\pi\over 2N+1}\right\}_{|q|\leq N},\,j=1,...,d\right\},\\
\end{array}
\end{equation}
and let $\nu$ be the normalized counting measure on $\Omega_N^d$:
$\nu(\{\omega\})=(2N+1)^{-d}$, $\omega\in\Omega_N^d$. Observe that
in view of {\sl R.2} the function $\delta(\cdot)$ maps
$\Omega_N^d$ into the unit disk $D=\{\zeta\in{\Bbb C}\mid\,
|\zeta|\leq 1\}$. Let $\mu$ be the distribution of values of
$\delta\big|_{\Omega_N^d}$, so that $\mu$ is the measure supported
by the finite set ${\cal M}=\{\zeta\mid\, \exists
\omega\in\Omega_N^d: \zeta=\delta(\omega)\}$, and
$\mu(\{\zeta\})=\sum\limits_{\omega\in\Omega_N^d:\delta(\omega)=\zeta}
\nu(\{\omega\}).$ Let also $
F(\alpha)=\mu\left(\left\{\zeta\mid\,\Re(\zeta)\geq
1-\alpha\right\}\right),\,\,\alpha\geq 0. $
\begin{lem}\label{lem:muandmu} There exists $c_2\in(0,1)$ such
that
\begin{eqnarray}\label{eq:calM}
\lefteqn{ {\cal M}\subset \widehat{\cal M}=\left\{\zeta\mid\,
|\zeta|\leq 1,|\Im(\zeta)|\leq
c_2^{-1}(1-\Re(\zeta))^{3/2}\right\},}\\
\lefteqn{F(\alpha)\leq
c_2^{-1}[\alpha^{d/2}+N^{-d}],\,\,0\leq\alpha\leq 2.}
\label{eq:massa}
\end{eqnarray}
\end{lem}
{\sl Proof.} (\ref{eq:calM}), (\ref{eq:massa}) are evident when
$\sum\limits_{\ell=1}^k \rho_\ell <1$, since then
$|\delta(\omega)|\leq 1-c_2$ for properly chosen $c_2$ and all
$\omega$. Thus, in the sequel we focus on the case of
$\sum\limits_{\ell=1}^k\rho_\ell=1$ (recall that
$\sum\limits_{\ell=1}^k\rho_\ell\leq1$ by {\sl R.2}). \\
{\sl 2$^0$.1)} Let ${\cal K}=\{\omega\in[-\pi,\pi]^d:
\delta(\omega)=1\}$. Since $\rho_\ell>0$,
$\sum\limits_\ell\rho_\ell=1$ and
$\delta(\omega)=\sum\limits_\ell\rho_\ell
\exp\{i\phi_\ell+\omega^T\alpha(\ell)\}$, a point $\omega\in{\cal
K}$ must satisfy the equations
\begin{equation}\label{eq:allequal}
\exp\{i[\phi_\ell+\omega^T\alpha(\ell)]\}=1\quad\forall (1\leq
\ell\leq k),
\end{equation}
whence $ \phi_\ell+\omega^T\alpha(\ell)\in 2\pi\bZ\quad \forall
(1\leq \ell\leq k). $ Since $\rank\{\alpha(\ell):1\leq\ell\leq
k\}= d$, the latter system of equations implies that ${\cal K}$
belongs to a set of the form $r+A\bZd$ with certain $d\times d$
nonsingular matrix $A$ (depending solely on ${\cal D}$). The
cardinality of the intersection of latter set with the cube
$[-\pi,\pi]^d$ does not exceed certain $c_3$. Thus, $\Card{\cal K}
\leq c_3$. \\
{\sl 2$^0$.2)} Let $\omega\in{\cal K}$, and let $d\omega\in{\Bbb
R}^n$ be such that $|d\omega|\leq 1$. Then
$$
\begin{array}{l}
\delta(\omega+d\omega)=\sum\limits_{\ell=1}^k \rho_\ell
\exp\{i[\phi_\ell+\omega^T\alpha(\ell)]\}\exp\{i(d\omega)^T\alpha(\ell)\}\\
\underbrace{=}_{a}\sum\limits_{\ell=1}^k
\rho_\ell \exp\{i(d\omega)^T\alpha(\ell)\}\\
\Rightarrow
|\delta(\omega+d\omega)|=\left|\sum\limits_{\ell=1}^k\rho_\ell
\exp\{i(d\omega)^T\alpha(\ell)\}\right|\\
\underbrace{\leq}_{b}\left|\sum\limits_{\ell=1}^k\rho_\ell\left(1+i(d\omega)^T\alpha(\ell)-{1\over
2}\left((d\omega)^T\alpha(\ell)\right)^2\right)\right|+c_4|d\omega|^3\\
\underbrace{=}_{b}\left|\sum\limits_{\ell=1}^k\rho_\ell\left(1-{1\over
2}\left((d\omega)^T\alpha(\ell)\right)^2\right)\right|+c_4|d\omega|^3
\underbrace{\leq}_{c}1-c_5|d\omega|^2+c_4|d\omega|^3
\end{array}
$$
(for $a$, see (\ref{eq:allequal}), $b$ is by
(\ref{eq:harmonic1}.$b$), $c$ is due to
$\rank\left(\{\alpha(\ell)\}_{\ell}\right)=d$). It follows that
with properly chosen $c_6$ one has
\begin{equation}\label{eq:sqrt}
\forall(\omega\in[-\pi,\pi]^d,|\delta(\omega)-1|\leq\alpha)\quad
\exists \bar{\omega}\in {\cal K}: |\omega-\bar{\omega}|\leq
c_6^{-1}\sqrt{\alpha}.
\end{equation} Since $\Card({\cal K})\leq c_3$ by 2$^0$.1) and
$|\delta(\omega)|\leq1$ for all $\omega$, we conclude that
\begin{equation}\label{eq:nunu}
\nu\left(\{\omega\in\Omega_N^d: |\delta(\omega)-1|\leq
\alpha\}\right)\leq c_7[\alpha^{d/2}+N^{-d}]\quad \forall
\alpha\leq 2.
\end{equation}
{\sl 2$^0$.3)} Now we can complete the proof of (\ref{eq:calM}),
(\ref{eq:massa}). Let $\bar{\omega}\in{\cal K}$, $d\omega\in{\Bbb
R}^d$, $|d\omega|\leq 1$. We have
\begin{equation}\label{eq:123456}
\begin{array}{l}
\delta(\bar{\omega}+d\omega)  =  \sum\limits_{\ell=1}^k \rho_\ell
\exp\{i[\phi_\ell+\bar{\omega}^T\alpha(\ell)]\}\exp\{i(d\omega)^T\alpha(\ell)\}\\
\underbrace{=}_{a}  \sum\limits_{\ell=1}^k \rho_\ell
\exp\{i(d\omega)^T\alpha(\ell)\}
  =  \sum\limits_{\ell=1}^k\rho_\ell\bigg(1+i(d\omega)^T\alpha(\ell)\\
  -{1\over
2}\left((d\omega)^T\alpha(\ell)\right)^2 -{i\over
6}\left((d\omega)^T\alpha(\ell)\right)^3+r_\ell(\omega,d\omega)\bigg),\\
\multicolumn{1}{r}{[|r_\ell(\omega,d\omega)|\leq c_{10}|d\omega|^4]}\\
\underbrace{=}_{b}  \sum\limits_{\ell=1}^k\rho_\ell\left(1-{1\over
2}\left((d\omega)^T\alpha(\ell)\right)^2-{i\over
6}\left((d\omega)^T\alpha(\ell)\right)^3+r_\ell(\omega,d\omega)\right)\\
\end{array}
\end{equation}
(for $a$, see (\ref{eq:allequal}), for $b$, see
(\ref{eq:harmonic1})). Taking into account that $\sum\limits_\ell
\rho_\ell=1$ and $c_{11}|d\omega|^2\leq \sum\limits_\ell \rho_\ell
\left((d\omega)^T\alpha(\ell)\right)^2\leq c_{12}|d\omega|^2$, we
conclude from (\ref{eq:sqrt}) combined with (\ref{eq:123456}) that
for properly chosen $c_{13}$ one has
$$
\omega\in[-\pi,\pi]^d\Rightarrow |\Im(\delta(\omega))|\leq
c_{13}(1-\Re(\delta(\omega)))^{3/2},
$$
and (\ref{eq:calM}) follows.  By (\ref{eq:calM}) one has $
|1-\delta(\omega)|\leq c_{14}(1-\Re(\delta(\omega))), $ so that
(\ref{eq:massa}) follows from (\ref{eq:nunu}).  \qed \\
{\sl 3$^0$.} Let $n$ be a positive integer, and let $T_n(\zeta)$
be the Tschebyshev polynomial of degree $n$. Recall that this
polynomial is defined as follows:
\begin{equation}\label{eq:Tscheb}
T_n(\zeta)={w^n+w^{-n}\over 2},\hbox{\ where
$w=\zeta+i\sqrt{1-\zeta^2}$.}
\end{equation}
In (\ref{eq:Tscheb}), the choice of the branch of  $\sqrt{\cdot}$
affects the value of $w$, but does not affect the value of
$w^n+w^{-n}$; since we intend to work with $\zeta$ from the unit
disk, so that $\Re(1-\zeta^2)>0$, in the calculations to follow we
deal with the main branch of $\sqrt{\cdot}$ in the closed right
half-plane. On the segment $[-1,1]$ of the real axis one has
$T_n(\zeta)=\cos(n\,\hbox{\rm acos}(\zeta))$, whence $T_n(1)=1$,
$T_n^\prime(1)=n^2$. From these relations it follows that the
function $ P_n(\zeta)={1-T_n(\zeta)\over n^2(1-\zeta)} $ is a
polynomial of degree $n-1$, and $ P_n(1)=1. $
\begin{lem}\label{lem:Pn} One has
\begin{equation}\label{eq:Pn}
\begin{array}{l}
p_n(\alpha)\equiv
\max\limits_\zeta\{|P_n(\zeta)|:\zeta\in\widehat{\cal
M},\Re(\zeta)=1-\alpha\}\\
\leq q_n(\alpha)=\cases{c_{15},&$0\leq\alpha\leq {1\over n^2}$\cr
{c_{15}(1+c_{15}\alpha)^n\over n^2\alpha},&${1\over
n^2}\leq\alpha\leq 2$\cr}.\\
\end{array}
\end{equation}
\end{lem}
{\sl Proof.} Let $\zeta=1-\alpha+i\beta\in\widehat{\cal M}$,  so
that
\begin{equation}
\label{eq:ma1} |\beta|\leq c_{16}\alpha^{3/2}.
\end{equation}
We have \begin{equation}\label{eq:ma2}
\begin{array}{l}
w\equiv\zeta+i\sqrt{1-\zeta^2}=1-\alpha+i\beta+i\sqrt{2\alpha-\alpha^2-2i(1-\alpha)\beta+\beta^2}\\
=1-\alpha+i\beta+i\sqrt{2\alpha}\sqrt{1-0.5\alpha+[0.5\beta-i(1-\alpha)](\beta/\alpha)}\\
=1+i\sqrt{2\alpha}+ r_1(\zeta),\quad |r_1(\zeta)|\leq
c_{17}\alpha\\
\end{array}
\end{equation}
(since $|\beta/\alpha|\leq c_{16}\sqrt{\alpha}$ by
(\ref{eq:ma2})). Note that completely similar considerations
demonstrate that
\begin{equation}\label{eq:ma1prime}
w^{-1}=\zeta-i\sqrt{1-\zeta^2}=1-i\sqrt{2\alpha}+r_2(\zeta),\quad
|r_2(\zeta)|\leq c_{17}\alpha.
\end{equation}
{\sl 3$^0$.1)} Assume, first, that $0\leq\alpha\leq{1\over n^2}$.
In this case from (\ref{eq:ma2}) it follows that $|1-w|\leq
\sqrt{2}n^{-1}$, whence, taking into account (\ref{eq:ma2}),
$$
\begin{array}{l}
|w^n-(1+n(w-1)+{n(n-1)\over 2}(w-1)^2)|\leq c_{17}(n|w-1|)^3
\leq c_{18}n^3\alpha^{3/2},\\
|w^{-n}-(1-n(w-1)+{n(n+1)\over 2}(w-1)^2)|\leq c_{17}(n|w-1|)^3
\leq c_{18}n^3\alpha^{3/2}\\
\Rightarrow \left|{w^n+w^{-n}\over2}-1\right|\leq
{n^2\over2}|w-1|^2+c_{18}n^3\alpha^{3/2} \leq
c_{19}(n^2\alpha+n^3\alpha^{3/2})
\leq c_{20}n^2\alpha.\\
\end{array}
$$
Thus,  one has
$|P_n(\zeta)|={\left|{w^n+w^{-n}\over2}-1\right|\over
n^2|\alpha-i\beta|}\leq c_{15}$,  as required in (\ref{eq:Pn}) for
the case of $0\leq\alpha\leq {1\over n^2}$. \\
{\sl 3$^0$.2)} Now consider the case of ${1\over
n^2}\leq\alpha\leq 2$. From (\ref{eq:ma2}), (\ref{eq:ma1prime}) it
follows that $|w|\leq 1+c_{21}\alpha$, $|w^{-1}|\leq
1+c_{21}\alpha$, whence $
|P_n(\zeta)|={\left|{w^n+w^{-n}\over2}-1\right|\over
n^2|\alpha-i\beta|}\leq {c_{22}(1+c_{21}\alpha)^n\over n^2\alpha},
$ as required in (\ref{eq:Pn}). \qed
\\
{\sl 4$^0$.} Let $Q(\zeta)={1+\zeta\over 2}$. It is immediately
seen that
\begin{equation}\label{eq:Qofzeta}
\zeta=1-\alpha+i\beta\in \widehat{\cal M}\Rightarrow
|Q(\zeta)|\leq 1-c_{23}\alpha \quad[c_{23}<{1\over 2}].
\end{equation}
Now let $c_{24}$ be a positive integer which is $\geq{c_{15}\over
c_{23}}$ (see (\ref{eq:Pn})). Consider the polynomial
$
S_n(\zeta)=P_n(\zeta)Q^{c_{24}n}(\zeta).
$
\begin{lem}\label{lem:Sn} For every positive integer $n$, the
polynomial $S_n(\zeta)$ possesses the following properties:
\begin{equation}\label{properties}
\begin{array}{l}
(a):\deg(S_n)\leq c_{25}n;\,
(b):S_n(1)=1;\\
(c):\max\limits_\zeta\{|S_n(\zeta)|:\zeta\in\widehat{\cal
M},\,\Re(\zeta)=1-\alpha\}\leq c_{15}\min\left[{1\over
n^2\alpha};1\right].\\
\end{array}
\end{equation}
\end{lem}
{\sl Proof.} Relations (\ref{properties}.$a-b$) are evident (take
into account that $P_n(1)=1$ and $\deg(P_n)\leq n$). To verify
(\ref{properties}.$c$), note that if
$\zeta=1-\alpha+i\beta\in\widehat{\cal M}$, then in view of
(\ref{eq:Pn}) one has
$$
\begin{array}{l}
0\leq\alpha\leq{1\over n^2}\Rightarrow
|S_n(\zeta)|\leq|P_n(\zeta)||Q(\zeta)|^{c_{24}n}\leq
|P_n(\zeta)|\leq c_{15};\\
{1\over n^2}\leq\alpha\leq2\Rightarrow
|S_n(\zeta)|\leq|P_n(\zeta)||Q(\zeta)|^{c_{24}n}
\underbrace{\leq}_{a} c_{15}{(1+c_{15}\alpha)^n\over
n^2\alpha}(1-c_{23}\alpha)^{c_{24}n}\\
\leq c_{15}{\exp\{c_{15}n\alpha\}\over
n^2\alpha}\exp\{-c_{23}c_{24}n\alpha\}
\underbrace{\leq}_{b} {c_{15}\over n^2\alpha}\\
\end{array}
$$
(for $a$, see (\ref{eq:Qofzeta}), $b$ is due to $c_{23}c_{24}\geq
c_{15}$). \qed\\
 {\sl 5$^0$.} Now we are ready
to complete the proof of Proposition \ref{prop:harmonic}. Given a
positive integer $n$, let us set
$
R_n(\zeta)=S_n^{d}(\zeta).
$
 In view of (\ref{properties}) one has
\begin{equation}\label{propertiesnew}
\begin{array}{l}
(a):\deg(R_n)\leq c_{26}n;\, (b):R_n(1)=1;\\
(c):\max\limits_\zeta\{|R_n(\zeta)|:\zeta\in\widehat{\cal
M},\,\Re(\zeta)=1-\alpha\}\leq r_n(\alpha)\\
\quad\quad\equiv c_{26}\min\left[{1\over
n^{2d}\alpha^d};1\right].\\
\end{array}
\end{equation}
Consider the filters $q^{(n)}(z)$ given by
$
q^{(n)}(\Delta)=R_n({\cal D}),\,\,n=0,1,...
$
By (\ref{propertiesnew}.$b$) and Lemma \ref{lem:harm1} we have
\begin{equation}\label{eq:result1}
\left.\begin{array}{r}
T\leq c_{27}N\\
1\leq n(T)\equiv \lfloor
c_{27}T\rfloor\\
(s)\in{\cal H}^t_N({\cal D})\\
\end{array}\right\}
\Rightarrow \left\{\begin{array}{l}\ord(q^{(n(T))})\leq T, \\
s_\tau=(q^{(n(T))}(\Delta)s)_\tau \, \forall(\tau:|\tau-t|\leq
c_{27}N).\\
\end{array}\right.
\end{equation}
By Parseval's equality, we have also (in what follows, $n=n(T)$)
\begin{equation}\label{eq:eqeqeq12}
\begin{array}{l}
|q^{(n)}|_2^2=\displaystyle{\int\limits_{\Omega^d_N}}
|R_n(\delta(\omega))|^2\nu(d\omega)
=\displaystyle{\int\limits_{\cal
M}}|R_n(\zeta)|^2\mu(d\zeta)\underbrace{\leq}_{a}
\displaystyle{\int\limits_0^2}\underbrace{r_n^2(\alpha)}_{\rho_n(\alpha)}dF(\alpha)\\
\end{array}
\end{equation}
with $a$ given by  (\ref{propertiesnew}.$c$), (\ref{eq:calM}) and
the definition of $F(\cdot)$. Let $\gamma$ be the measure on
$[0,2]$ defined by $
G(\alpha)\equiv\gamma([0,\alpha])=c_2^{-1}(\alpha^{d/2}+N^{-d}), $
so that \begin{equation}\label{eq:eqeqeq13} F(\alpha)\leq
G(\alpha)\equiv \gamma([0,\alpha])\quad \forall \alpha\in[0,2]
\end{equation}
 (see
(\ref{eq:massa})). We have
\begin{equation}\label{eq:eqeqeq14}
\begin{array}{l}
\displaystyle{\int\limits_0^2}\rho_n(\alpha)dF(\alpha) =
\rho_n(2)-
\displaystyle{\int\limits_0^2}\rho_n^\prime(\alpha)F(\alpha)d\alpha
 \underbrace{\leq}_{a} \rho_n(2)\\
 -\displaystyle{\int\limits_0^2}\rho_n^\prime(\alpha)G(\alpha)d\alpha
 =
 \rho_n(2)-\rho_n(2)G(2)+\displaystyle{\int\limits_0^2}\rho_n(\alpha)\gamma(d\alpha)\\
 \underbrace{\leq}_{b}
  \displaystyle{\int\limits_0^2}\rho_n(\alpha)\gamma(d\alpha)
\underbrace{=}_{c}
c_2^{-1}\bigg[c_{28}\displaystyle{\int\limits_0^2}{\min}^2\left[n^{-2d}\alpha^{-d},1\right]
\alpha^{{d\over2}-1}d\alpha\\
+\rho_n(0)N^{-d}\bigg] \underbrace{\leq}_{d}
c_{30}\left[N^{-d}+n^{-d}\right] \leq c_{31}(2T+1)^{-d}
\end{array}
\end{equation}
($a$ holds since $\rho_n(\cdot)$ is nonincreasing, see
(\ref{propertiesnew}.$c$), and by (\ref{eq:eqeqeq13}), $b$ holds
since $c_2\in(0,1)$, see Lemma \ref{lem:muandmu}, $c$ is by
(\ref{propertiesnew}.$c$) and (\ref{eq:eqeqeq12}), $d$ is due to $
n=n(T)=\lfloor c_{27}T\rfloor$). Combining (\ref{eq:eqeqeq12}) and
(\ref{eq:eqeqeq14}), we conclude that
\begin{equation}\label{eq:result2}
|q^{(n(T))}|_2\leq c_{32}(2T+1)^{-d/2}.
\end{equation}
From (\ref{eq:result1}) and (\ref{eq:result2}) we conclude that if
$L=\lfloor c_{33}N\rfloor$ and $T\leq L$ is such that $n(T)\equiv
\lfloor c_{27}T\rfloor\geq1$, then
\begin{equation}\label{eq:finalres}
\exists q^{(T)}\in C_T(\bZd): \left\{\begin{array}{l}
|q^{(T)}|_2\leq c_{32}(2T+1)^{-d/2},\\
s_\tau=(q^{(T)}(\Delta)s)_\tau \,\forall(\tau,|\tau-t|\leq
L,(s)\in{\cal H}^t_N({\cal D}))\\
\end{array}\right. \end{equation} (indeed, one can choose, as a required $q^{(T)}$,
the filter  $q^{(n(T))}$). Setting $q^{(T)}(z)\equiv 1$ for
$T<{1\over c_{27}}$, we enforce the validity of
(\ref{eq:finalres}) for all $T$, $0\leq T\leq L$. Thus, ${\cal
H}^t_N({\cal D})\subset {\cal F}^t_{\lfloor
c_{29}L\rfloor}(0,c_{34})$. \qed \subsubsection{Proof of
Proposition \protect{\ref{propharm}}}
\begin{lem}\label{lem:harm} Let $f\in {\cal H}^+(M)$ be a deterministic function, let $N\leq M/2$, and let
$t\in \bZ^d$, $|t|\leq N$. Consider the ``discrete box''
$B^t_N=\{\tau\in\bZ^d:|\tau-t|\leq N\}$, and let $\phi$ be a
deterministic function on $B^t_N$ which coincides with $f$ on the
``discrete boundary'' $\partial
B^t_N\equiv\{\tau\in\bZ^d:|\tau-t|=N\}$ of $B^t_N$ and is
``discrete harmonic'':
$
\tau\in\bZ^d, |\tau-t|<N\Rightarrow \phi_\tau={1\over
2d}\sum\limits_{{\epsilon=(\epsilon_1,...,\epsilon_d)\atop
|\epsilon_1|=...=|\epsilon_d|=1}}\phi_{\tau+\epsilon}.
$
Then
\begin{equation}\label{harmbound}
\tau\in B^t_N \Rightarrow |f(\tau)-\phi_\tau|\leq
c_1\|f\|_{\infty,2M}N^{-2}
\end{equation}
(from now on, $c_i$ are positive absolute constants).
\end{lem}
{\sl Proof.} First, we should prove that the ``discrete harmonic''
function $\phi$ on $B^t_N$ which coincides with $f$ on $\partial
B^t_N$ does exist. This fact is well known; we present here its
proof just for the sake of completeness. Let $\psi$ be a function
on $\partial B^t_N$. Consider the following random walk on
$B^t_N$: arriving for the first time at a point $\tau$ from
$\partial B^t_N$, we pay penalty $\psi(\tau)$ and terminate; from
an ``interior point'' $\tau\in \inter B^t_N\equiv B^t_N\backslash
\partial B^t_N$ we make a random step of length 1 along one of the coordinate
axes, choosing every one of $2d$ possible steps with probability
$1/(2d)$. It is immediately seen that the expected penalty payed
at the termination, treated as a function of the initial state, is
a discrete harmonic function with the boundary values $\psi$.
\\
Now, since $|t|\leq N$ and $2N\leq M$, the function $f$ is
harmonic in the ``continuous box'' $D^t_{2N}=\{\tau\in{\Bbb
R}^d:|\tau-t|\leq 2N\}$, and the uniform norm of $f$ in this
square does not exceed $\|f\|_{\infty,2M}$.  From the standard
results on harmonic functions it follows that
\begin{equation}\label{asnem115}
\forall (\tau\in D^t_N): \left|{\partial^\kappa\over\partial
x_j^\kappa}f(\tau)\right|\leq
c_{2}\|f\|_{\infty,2M}N^{-\kappa},\,\kappa=1,2,3,4,\,j=1,...,d.
\end{equation}
Consequently, for the basic orths $e_j$, $j=1,...,d$ we have
$$
\tau\in D^t_N,|s|\leq 1\Rightarrow \left|f(\tau+se_j)-
\sum\limits_{\kappa=0}^3
{1\over\kappa!}{\partial^\kappa\over\partial
x_j^\kappa}f(\tau)s^\kappa\right|\leq
c_3|s|^4\|f\|_{\infty,2M}N^{-4}.
$$
Since $f$ is harmonic, we conclude that
\begin{equation}\label{asnem1200}
|({\cal D}f)_\tau|\leq c_4\|f\|_{\infty,2M}N^{-4}, \,\tau\in
B^t_N.
\end{equation}
Now let $h=f\big|_{\bZ^d}-\phi\in C(B^t_N)$ and let
$
h^\pm_\tau=h_\tau\pm{2c_4\|f\|_{\infty,2M}\over
N^4}\sum\limits_{j=1}^d(\tau_j-t_j)^2.
$
Taking into account (\ref{asnem1200}) and the fact that $\phi$ is
discrete harmonic, we have for $\tau\in\inter B^t_N$:
$$
({\cal D}h^+)_\tau=({\cal D}h)_\tau + {2c_4\|f\|_{\infty,2M}\over
N^4}>0,\, ({\cal D}h^-)_\tau=({\cal D}h)_\tau -
{2c_4\|f\|_{\infty,2M}\over N^4}<0,
$$
whence both the maximum of $h^+$ and the minimum of $h^-$ over
$B^t_N$ are attained at $\partial B^t_N$. Since at the discrete
boundary of $B^t_N$ we have $f=\phi$ and therefore $h^+\leq
4c_4\|f\|_{\infty,2M}N^{-2}$, we conclude that $ \tau\in B^t_N
\Rightarrow h_\tau\leq h^+_\tau\leq \max\limits_{\tau\in\partial
B^t_N} h^+_\tau \leq 2dc_4\|f\|_{\infty,2M}N^{-2}. $ By similar
reasons, $ \tau\in B^t_N \Rightarrow h_\tau\geq h^-_\tau\geq
\min\limits_{\tau\in\partial B^t_N} h^-_\tau \geq
-2dc_4\|f\|_{\infty,2M}N^{-2}. $ \qed \\
 Now let $|t|\leq M/8$ and
$L\leq M/8$. Given $T$, $0\leq T\leq L$, and applying Proposition
\ref{prop:harmonic}, we can build filter $q^{(T)}\in C_T(\bZ^d)$
such that
\begin{equation}\label{asnem125}
|q^{(T)}|_2\leq c_5(2T+1)^{-1},\,\, \phi_\tau =
\sum\limits_{|\nu|\leq T} \phi_{\tau-\nu}q^{(T)}_\nu \,\,\forall
(\tau:|\tau-t|\leq L)
\end{equation}
for every $\phi$ which is discrete harmonic in the discrete box
$B^t_{2L}$. Now let $f\in {\cal H}(M,R)$. Applying Lemma
\ref{lem:harm}, we can find function $\phi$ which is discrete
harmonic in the box $B^t_{2L}$ and such that
$|\phi_\tau-f_\tau|^2\leq c_6^2\|f\|_{\infty,2M}^2L^{-4}$ for
$\tau\in B^t_{2L}$. From (\ref{asnem125}) it now follows that
$$
\begin{array}{l}
\forall (\tau:|\tau-t|\leq L):
\left[E\left\{|f_\tau-\sum\limits_{|\nu|\leq T}
f_{\tau-\nu}q^{(T)}_\nu|^2\right\}\right]^{1/2}\\
\leq c_6\underbrace{[E\{\|f\|_{\infty,2M}^2\}]^{1/2}}_{\leq
R}L^{-2}(1+|q^{(T)}|_1) \leq
c_6RL^{-2}(1+|q^{(T)}|_2(2T+1)^{d/2})\\ \leq c_8RL^{-2}
\leq c_9R(2T+1)^{-d/2}\\
\end{array}
$$
(recall that $d\leq 4$). \qed
\par
{\sl The proof of Proposition \ref{prophomogpred}} is completely similar to that of Proposition
\ref{prophomog}.

\end{document}